\documentclass[times,3p]{elsarticle}

\usepackage{amssymb}
\usepackage{latexsym}
\usepackage{amsmath}

\usepackage{ifthen}
\usepackage{ulem}

\usepackage{lipsum}
\usepackage{amsfonts}
\usepackage{graphicx}
\usepackage{epstopdf}

\usepackage{multirow}
\usepackage{subcaption}
\usepackage{hyperref}
\usepackage{cleveref}

\usepackage{tikz}
\usetikzlibrary{matrix}
\usetikzlibrary{shapes,patterns,snakes,arrows}
\usetikzlibrary{arrows.meta}
\usetikzlibrary{positioning,3d,fit,calc}
\usetikzlibrary{graphs}

\usepackage{todonotes}

\usepackage{xcolor}
\definecolor{darkred}{rgb}{.9,0,.05}

\usepackage{algorithmic}
\ifpdf
\DeclareGraphicsExtensions{.eps,.pdf,.png,.jpg}
\else
\DeclareGraphicsExtensions{.eps}
\fi

\newcommand{\mycomment}[1]{%
\ifthenelse{\isodd{\value{page}}}{%
\normalmarginpar%
\marginpar{\tiny {#1}}%
}{%
\reversemarginpar%
\marginpar{\tiny {#1}}%
}}%

\usepackage[symbol]{footmisc}

\newcommand{\takeout}[1]{}

\begin{document}
	

\begin{frontmatter}
	\title{Learning the solution operator of two-dimensional incompressible Navier-Stokes equations using physics-aware convolutional neural networks\footnote[7]{This work is based on the doctoral dissertation of the first author at the Faculty of Mathematics and Natural Science of the University of Cologne.}}
	
	\author[1]{Viktor Grimm \corref{cor1}}
	\cortext[cor1]{Email addresses: viktor.grimm@uni-koeln.de, a.heinlein@tudelft.nl, axel.klawonn@uni-koeln.de}
	
	\author[3]{Alexander Heinlein \corref{cor1}} 
	
	\author[1,2]{Axel Klawonn \corref{cor1}}
	
	\address[1]{Department of Mathematics and Computer Science, University of Cologne, Weyertal 86-90, 50931 K\"oln, Germany}
	\address[2]{Center for Data and Simulation Science, University of Cologne, 50923 K\"oln, Germany}
	\address[3]{Delft Institute of Applied Mathematics, Delft University of Technology, Mekelweg 4, 2628 CD Delft, Netherlands.}
	
	
	\begin{abstract}
		In recent years, the concept of introducing physics to machine learning has become widely popular. 
		Most physics-inclusive ML-techniques however are still limited to a single geometry or a set of parametrizable geometries.
		Thus, there remains the need to train a new model for a new geometry, even if it is only slightly modified.
		With this work we introduce a technique with which it is possible to learn approximate solutions to the steady-state Navier--Stokes equations in varying geometries without the need of parametrization. 
		This technique is based on a combination of a U-Net-like CNN and well established discretization methods from the field of the finite difference method.
		The results of our physics-aware CNN are compared to a state-of-the-art data-based approach. 
		Additionally, it is also shown how our approach performs when combined with the data-based approach.
	\end{abstract}
	
\end{frontmatter}

\section{Introduction}
\label{sec:introduction}
Fluid behavior is important in various fields such as civil, mechanical, and biomedical engineering, aerospace, meteorology, and geosciences. 
The governing equations for fluid behavior are typically the Navier-Stokes equations, which are solved using discretization approaches like finite difference, finite volume, or finite element methods. 
However, such computational fluid dynamics (CFD) simulations can be computationally intensive, especially for turbulent flow and complex geometries, and changing the geometry requires recomputing the entire simulation. 
Hence, there is a need for a quick surrogate model for CFD simulations. 
Such surrogate models encompass a variety of approaches, including linear reduced order models~\cite{frs_liii_1901, LMQR:2014:ROMNS}, such as reduced basis~\cite{quarteroni_reduced_2016} and proper orthogonal decomposition~\cite{rathinam_new_2003} models, as well as neural network-based models~\cite{fresca_comprehensive_2021}, like convolutional neural networks (CNNs)~\cite{eichinger_surrogate_2022, franco_approximation_2023, GLI:2016:NFA, lee_model_2019, maulik_reduced-order_2021} and neural operators~\cite{kovachki_neural_2022, lu_deeponet_2021}.

In present work, we focus on using neural networks as an approximation for CFD simulations.
Instead of relying on a large dataset, we leverage the known governing equations of fluids to construct a physics-aware loss function and train our model to satisfy these equations discretely.
This approach has recently become increasingly popular and was applied to dense neural networks (DNN) to solve partial differential equations (PDEs) with little training data \cite{RPK:2019:MID} or without training data \cite{SGPW:2020:SMF} as well as inverse problems with limited training data \cite{GHKLW:2020:CORONA, RPK:2019:MID}.
More recently, this idea was also applied to convolutional neural networks (CNN) by using physics-aware loss functions to solve PDEs  \cite{BSDM:2019:StackPCNN, F:2021:PINNCNN, GSW:2021:PGN, RRL:2021:PhyCRNet, SFG:2018:WSL, ZLS:2020:PhyCNN}, upscale and denoise solutions \cite{GSW:2021:SuperResPGN, KRM:2022:PICNNSuperRes}, generally improve the predictive quality of a model \cite{SH:2019:PCNNFaultDiag, WYX:2018:FPCNN}, or learn PDEs from data \cite{long_pde-net_2019, long_pde-net_2018}.
For a comprehensive overview on scientific machine learning (SciML), we refer to \cite{PRR:2022:PINNsOverview, WJXSK:2021:ReviewSM}.

However, previous physics-informed machine learning approaches have imposed geometric constraints, such as rectangularity or parametrizability, or have been limited to specific geometries or even to a single geometry.
In practice, these conditions are usually not met. 
Furthermore, the exact geometry is often unknown or at least not known in sufficient detail.
This is, for example, the case with medical imaging procedures.
Therefore, we aim to develop a CNN that is capable of learning the flow field and pressure under physics constraints without labeled data for, with some restrictions, arbitrary geometries.
We explicitly do not use methods such as coordinate transformations, which would limit us to specific geometries.
This approach aims to advance existing methods for more realistic applications, and to the authors' knowledge, it is the first attempt to use a single CNN for multiple irregular geometries without labeled data.

The rest of this paper is organized as follows.
We first define our stationary boundary value problem in \cref{sec:flow-problem}. 
Then, we introduce CNNs as surrogate models in \cref{sec:cnn} and directly afterwards extend this framework to physics-aware label-free learning of PDE solutions in \cref{sec:method} and its application to the Navier-Stokes equations in \cref{sec:ns}.
Next, in \cref{sec:model}, we explain the architecture of our surrogate model.
The creation of the training data used here is described in \cref{sec:trainingdata}. 
We present results for the data-based and physics-aware approach in  \cref{sec:results}. 
Finally, we draw conclusions in \cref{sec:conclusion}.

\subsection{Physics--Informed Machine Learning} \label{sec:PIML}
The method described in this paper can be categorized as physics-informed Machine Learning (ML)~\cite{meng2022physics}.
This term is used to describe ML methods for which prior knowledge, often known physical laws, are used for training.
Since we use the governing equations to train our physics-aware CNN, our method can also be referred to as a Physics-Informed Neural Network (PINN).
However, this would allow for possible confusion with the classical PINN approach, introduced in \cite{raissi2017physics, RPK:2019:MID} and pioneered in \cite{lagaris_artificial_1998}, where a physics-based loss function is constructed by clever use of automatic differentiation, classically for dense neural networks.
In particular, the classical PINN approach requires an NN to directly approximate the solution function, i.e., the discretization is done by the NN.
In contrast, in our approach, we use a finite difference-based discretization and predict the coefficients using a CNN.
Consequently, we approximate the residuals using finite differences.
Note, however, that we could also use finite element or finite volume techniques.

\section{Model problem} \label{sec:flow-problem}

\begin{figure}[t]
	\centering
	\begin{tikzpicture}[scale=1.5]
	\node (geometry_lower_left) at (0.0,0.0) {}; 
	\node (geometry_upper_left) at (0.0,3.0) {}; 
	\node (geometry_upper_right) at (6.0,3.0) {};  
	\node (geometry_lower_right) at (6.0,0.0) {}; 
	
	\draw[fill={rgb:black,1;white,2}, line width=0.5mm] (geometry_lower_left.center) rectangle (geometry_upper_right.center);
	
	\draw[darkred, line width=0.5mm] (geometry_upper_right.center) -- node[rotate=-90,anchor=south,thick] (outflow) {$\partial\Omega_{\text{out}}$} (geometry_lower_right.center);
	\draw[blue, line width=0.5mm] (geometry_upper_left.center) -- node[rotate=90,anchor=south,thick] (inflow) {$\partial\Omega_{\text{in}}$} (geometry_lower_left.center);
	
	\node (no-slip) at (1.25, 0.25) {$\partial\Omega_{\text{wall}}$};
	
	\draw[fill=white, thick] (3.30500756, 0.98558888) -- (3.01580547, 1.33185977) -- (3.87749112, 2.02937552) -- (4.33515942, 2.06316199) -- (3.30500756, 0.98558888);
	
	\draw[dashed, thick] (0.75,2.25) -- (5.25,2.25); 
	\draw[dashed, thick] (5.25,2.25) -- (5.25,0.75); 
	\draw[dashed, thick] (0.75,0.75) -- (5.25,0.75); 
	\draw[dashed, thick] (0.75,0.75) -- (0.75,2.25); 
	\node[thick] (P-text) at (3.5, 1.45) {$P$};
	
	\draw[<->] (3, 0.05) -- (3, 0.7) node(lower-distance-arrow)[midway] {};
	\node (lower-distance-text) [right=0.05cm of lower-distance-arrow] {$0.75$};
	
	\draw[<->] (5.3, 2.0) -- (5.95, 2.0) node(outflow-distance-arrow)[midway] {};
	\node (outflow-distance-text) [above=0.05cm of outflow-distance-arrow] {$0.75$};
	
	\draw[|-|, thick] (6.6, 0) -- (6.6, 3) node(height-arrow)[midway] {};
	\node (height-text) at (6.75, 1.5) {$3$};
	
	\draw[|-|, thick] (0, 3.2) -- node[midway,above] {$6$} (6, 3.2) node(width-arrow)[midway] {};
	\end{tikzpicture}
	\caption{Example of a channel geometry $\Omega$ with a star-shaped obstacle. The obstacle is confined within a box (dashed line) with a distance of $0.75$ to the boundary of $\Omega$.
	\label{fig:showcase_domain}
	}
\end{figure}
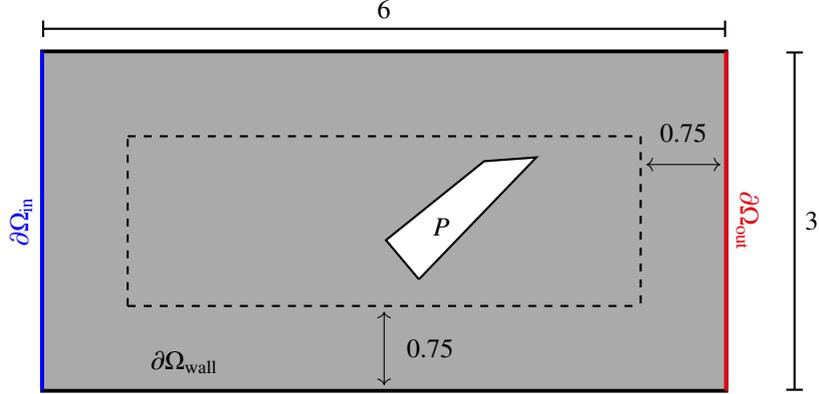

We consider the stationary Navier--Stokes equations describing incompressible Newtonian fluids with constant density
\begin{equation} \label{eq:ns}
	\begin{aligned}
		(\vec{u} \cdot \nabla) \vec{u} - \nu \Delta \vec{u} + \nabla p &= 0  \quad \text{in } \Omega, \\
		\nabla \vec{u} &= 0 \quad \text{in } \Omega,
	\end{aligned}
\end{equation}
where $\vec{u}$ is the velocity field, $p$ the pressure, and $\nu$ the kinematic viscosity. 

In our experiments, we consider rectangular channels $\Omega = [0,6] \times [0,3]$ from which we have cut a star-shaped obstacle $P$; cf.~\cref{fig:showcase_domain} for an example. This design is inspired by~\cite{GLI:2016:NFA,eichinger_surrogate_2022}. 
We apply the following boundary conditions: On the inlet $\partial\Omega_{in} := 0 \times [0,3]$, we prescribe a constant inflow velocity $u = \left( 3, 0 \right)^\top$, and at the outlet $\partial\Omega_{out} := 6 \times [0,3]$, we fix the pressure to $p = 0$. The lower and upper parts of the boundary $0 \times [0,6]$ and $3 \times [0,6]$, respectively, correspond to walls, and hence we enforce no-slip conditions $u = \left( 0, 0 \right)^\top$. Finally, we choose $\nu = 5 \cdot 10^{-2}$.
Depending on the shape and position of the obstacle, this setup leads to strongly varying flow patterns, making the channel problem a challenging benchmark for a CFD surrogate model.

\section{Surrogate models based on CNNs} \label{sec:cnn}

Let us first discuss the approach to construct surrogate models via CNNs from~\cite{GLI:2016:NFA,eichinger_surrogate_2022}, which is the basis for this work. CNNs~\cite{lecun1989generalization} are  artificial neural networks (ANNs) that employ linear transformations based on discrete convolutions within the network layers
making them well suited for structured temporal or spatial data, where neighboring coefficients correspond to neighboring points in space or time, respectively. 

CNNs are therefore also suitable for the approximating the solutions of partial differential equations on a structured tensor product grid: even though interaction is typically global, it is strongest for neighboring nodes. If the data structure is based on unstructured grid, graph convolutional networks (GCNs)~\cite{kipf2016semi} can be employed as an alternative.
Here, we only consider tensor product-structured data, and therefore, restrict ourselves to classical CNNs.

In the approach from~\cite{GLI:2016:NFA,eichinger_surrogate_2022}, a CNN that maps from the geometry of the computational domain to solution field(s) of the corresponding boundary value problem is trained; here, we specifically aim at predicting the velocity and pressure fields satisfying the Navier--Stokes equations~\cref{eq:ns}. In order to be able to employ standard CNNs, the geometry and solution fields are therefore interpolated to a tensor product grid. In two dimensions, the resulting data has a simple matrix structure; see~\cref{fig:showcase_data} for an exemplary pair of input and output data.
As can be seen, due to the matrix structure of the input and output data, they can be directly identified as pixel images. This also allows us to use a large variety of techniques from the application of CNN models to image data.

Let us now give a formal introduction of the approach. 
Therefore, let $I_g \in \mathbb{R}^{w \times h}$ be the pixel image matrix representing some computational domain $\Omega_g$; $g$ indicates a generic index for a specific geometry. 
Moreover, let $u_g^i \in \mathbb{R}^{w \times h}$ be a matrix representation of the $i$th component of the solution field of the boundary value problem solved on the computational domain $\Omega_g$.
Here, $w$ and $h$ correspond to the width and height of the pixel images, as well as the number of interpolation nodes in the $x$ and $y$ directions.
By assembling the tensors $u_g^i$, $i=1,\ldots,d$, we obtain a third-order tensor $u_g \in \mathbb{R}^{d \times w \times h}$, assuming that all solution components are defined on the same pixel grid.
For a simple diffusion equation $d$ is one, and for the Navier--Stokes equations in two dimensions, we have two velocity components and one pressure component, such that $d$ is three. 
In analogy to pixel images, each component of the solution field is regarded as one channel of the output image.

Our goal is to train a CNN that approximates the solution operator
$$
\begin{aligned}
	\mathcal{U} : \mathbb{R}^{w \times h} & \to \mathbb{R}^{d \times w \times h} \\
	I_g & \to u_g,
\end{aligned}	
$$
that is, the operator that maps a pixel representation of the geometry of the computational domain to a pixel representation of the solution of the corresponding boundary value problem. 
Hence, our approach can be seen as an example of operator learning; cf.~the related DeepONet~\cite{lu_deeponet_2021} and Fourier neural operator~\cite{li_fourier_2021} approaches, which employ different network architectures.

Let us denote the CNN model by $f_{N\!N}^\Psi$, where $\Psi$ are the trainable network parameters. 
In~\cite{GLI:2016:NFA,eichinger_surrogate_2022}, a CNN $f_{N\!N}$ has then been trained to approximate the solution operator $\mathcal{U}$ in a purely data-driven way. 
In particular, high-fidelity simulation data has been employed as the reference data $u_g$, and the model has been trained to minimize the mean squared error (MSE) between the model output and $u_g$. 
This corresponds to the minimization problem:
\begin{equation} \label{eq:mse}
	\arg\min_{\Psi} \frac{1}{\left| T \right|} \sum_{g \in T} \left\| f_{N\!N}^\Psi \left( I_g \right) - u_g \right\|^2,
\end{equation}
where $T$ is a set of geometries used as training data. Note training the model $f_{N\!N}^\Psi$ with this loss function requires the availability of reference data $u_g$; this means that a large number of measurement or high fidelity simulation data has to be available before the model training. 

It remains to discuss how to construct $I_g$ and $u_g$ for a specific geometry $g$. The approach is not restricted to a specific image representation of the geometry $I_g$, and in~\cite{eichinger_surrogate_2022, GLI:2016:NFA}, a binary or signed distance function (SDF)-based image of the geometry have been employed. 
It can be observed that the SDF input yields slightly better results. 
However, it comes at a computational cost, and the computation of the exact SDF input image requires precise knowledge about the boundary of the geometry; in practical applications, for instance, when the geometry is only known from medical image data, the SDF function can only be computed approximately based on the available image data.
A binary input image can be generated more easily by checking if the center or most of the volume of each pixel lies within the computational domain $\Omega_g$.
Here, we only consider binary input images.

The output pixel images $u_g$ can be constructed from a reference solution $\hat u$ by an interpolation operator, for instance, point-wise interpolation in the center points of the pixels or by averaging over the pixels (Cl\'ement-type interpolation). Here, $\hat u$ could be high-fidelity simulation or measurement data. We discuss the data processing for this paper in more detail in~\cref{sec:trainingdata}.

Next, we introduce the main novelty of this paper, that is, our approach to replace the data-based loss function~\cref{eq:mse} by a physics-aware loss function that requires only the knowledge of the PDE.

\section{A physics-aware surrogate model based on CNNs}
\label{sec:method}
In this section, we extend the approach from Section \ref{sec:cnn} by incorporating knowledge of the mathematical PDE model during network training. 
While the loss function in \cref{eq:mse} relies on reference data, our novel approach only requires a mathematical formula for the PDE residual, although a combination of both loss functions is possible. 
By considering the network output as a discrete finite difference solution on the same grid, we can approximate the PDE residual using finite difference stencils. 
This method can be extended to other discretization approaches on a structured grid, such as finite element or finite volume methods, but we focus on finite difference discretization for simplicity.

In this section, we extend the approach in~\cref{sec:cnn} by including knowledge of the mathematical PDE model into network training. 
Whereas the loss function in~\cref{eq:mse} relies on reference data, our novel approach only requires a mathematical formula for the residual of the PDE; a combination of both loss functions is also possible. 
By considering the network output as a discrete finite difference solution on the pixel image grid, we can approximate the PDE residual using finite difference stencils.
This method can be extended to other discretization approaches on a structured grid, such as finite element or finite volume methods, but we focus on finite difference discretization for simplicity and leave other discretization approaches to future work.

Like the data-driven approach described in~\cref{sec:cnn}, our new approach yields a surrogate model capable of predicting solutions for a range of geometries of the computational domain. 
It can be trained without using reference data or in combination with reference data. 
To introduce the approach, we first explain how to apply it to a stationary diffusion equation; see also~\cite{GHK:2023:ShortNote} for preliminary results for the stationary diffusion equation. 
Then, we discuss specifically how to apply the approach to the two-dimensional Navier--Stokes equations and the handling of boundary conditions.

\subsection{Finite differences and discrete convolutions}
\label{sec:FD_and_convolutions}
In order to derive the implementation of finite difference stencils based on the discrete convolution respectively cross-correlation operation, which is generally used in convolutional neural networks, we first consider a simple stationary diffusion problem: find the function $u$ such that
\begin{equation} \label{eq:pde}
	\begin{aligned}
		\frac{\partial^2 u}{\partial x^2}
		+
		\frac{\partial^2 u}{\partial y^2} 
		& = f \quad \text{in } \Omega = [0,1]^2,
	\end{aligned}
\end{equation}
where we assume that $f$ is a continuous function on $\Omega$. 
Later, in~\cref{sec:ns}, we discuss the application to the Navier--Stokes equations, which are our focus in this work.

Now, let us introduce a uniform grid $\Omega_h = \left\lbrace \boldsymbol{x}_{ij} \mid 1 \leq i,j \leq n+1 \right\rbrace$, with $\boldsymbol{x}_{ij} =  \left(\left(i-1\right) h, \left(j-1\right)  h\right)$ and $h = \frac{1}{n}$, and let
$
	U^h
	=
	\left(U_{ij}^h\right)_{ij}
	\in
	\mathbb{R}^{\left(n+1\right) \times \left(n+1\right)}
$
with $U_{ij}^h \approx u \left(\boldsymbol{x}_{ij}\right)$ the matrix representation of a discretization of $u$. Then,
$
	u^h
	=
	\left(u_{i}^h\right)_{i}
	\in
	\mathbb{R}^{\left(n+1\right) \cdot \left(n+1\right)}
$
with 
\begin{equation} \label{eq:discrete_u}
	u_{\left(j-1\right) \cdot \left(n+1\right) + i}^h = U_{ij}^h
\end{equation}
is the corresponding vector representation in lexicographical order. Then, discretizing~\cref{eq:pde} using central differences involves the approximation
\begin{equation} \label{eq:fdm_cd}
	\frac{\partial^2 u}{\partial x^2} 
	\approx
	\frac{U_{i+1,j}^h - 2U_{i,j}^h + U_{i-1,j}^h}{h^2}
	\quad \text{and} \quad
	\frac{\partial^2 u}{\partial y^2} 
	\approx
	\frac{U_{i,j+1}^h - 2U_{i,j}^h + U_{i,j-1}^h}{h^2},
\end{equation}
which could also be rewritten in terms of the entries of $u^h$ using~\cref{eq:discrete_u}.	This leads to a  linear system of equations
\begin{equation}
	\label{eq:laplace_lgs}
	A u^h = f^h\text{.}
\end{equation}
For the right-hand side vector, let
$
	F^h
	=
	\left(F_{ij}^h\right)_{ij}
	\in
	\mathbb{R}^{\left(n+1\right) \times \left(n+1\right)}
$
be the matrix representation with $F_{ij}^h = f \left(\boldsymbol{x}_{ij}\right)$ and
$
	f^h
	=
	\left(f_{i}^h\right)_{i}
	\in
	\mathbb{R}^{\left(n+1\right) \cdot \left(n+1\right)}
$
be the corresponding vector representation, where 
\begin{equation} \label{eq:discrete_f}
	f_{\left(j-1\right) \cdot \left(n+1\right) + i}^h = F_{ij}^h.
\end{equation}

We can observe that 
\begin{equation}
	\label{eq:laplace_lgs_equal_to_conv}
	A u^h = f^h 
	\quad
	\Leftrightarrow
	\quad
	U^h \ast K = F^h,
\end{equation}
where $\ast$ is the cross-correlation operation and
\begin{equation} \label{eq:fdm:stencil}
	K 
	=
	\frac{1}{h^2}
	\begin{pmatrix}
		K_{-1, -1} & K_{-1, 0} & K_{-1, 1} \\
		K_{0, -1} & K_{0, 0} & K_{0, 1} \\
		K_{1, -1} & K_{1, 0} & K_{1, 1}
	\end{pmatrix}
	=
	\frac{1}{h^2}
	\begin{pmatrix}
		0 & 1 & 0 \\
		1 & -4 & 1 \\
		0 & 1 & 0
	\end{pmatrix},
\end{equation}
is the kernel which corresponds to the finite difference stencil of the central difference scheme~\cref{eq:fdm_cd}. The cross-correlation is the linear transformation implemented in the convolutional layers in CNNs of current state-of-the-art machine learning libraries; cf.~\cite[section~9.1]{goodfellow2016deep}. 

In general form, the cross-corelation is given by
\begin{equation} \label{eq:cross-correlation}
	\left( I \ast K \right)_{ij}
	=
	\sum_{m} \sum_{n} I_{i+m,j+n} K_{m,n},
\end{equation}
where $I$ is some matrix and $K$ is, again, a kernel matrix, such as the one given in~\cref{eq:fdm:stencil}. As by convention, we omit the range of the sums, and regard each matrix coefficient as zero which is outside the range of indices. Note that flipping the kernel in~\cref{eq:cross-correlation} yields the discrete convolution 
$$	
	\left( I \tilde{\ast} K \right)_{ij}
	=
	\sum_{m} \sum_{n} I_{i-m,j-n} K_{m,n};
$$
see, for instance,~\cite[Section~9.1]{goodfellow2016deep}. Since, in a CNN, the entries of the kernel are generally trainable, the cross-correlation and the discrete convolution are equivalent in that sense.

By numbering the rows and columns of $K$ from $-1$ to $1$, as done in~\cref{eq:fdm:stencil}, 
we can easily show the equality of the left hand sides in~\cref{eq:laplace_lgs_equal_to_conv}:
\begin{align*}
	\left( A u^h\right)_{\left(j-1\right)(n+1)+i} 
	\stackrel{\cref{eq:fdm_cd}}{=} & ~\frac{1}{h^2} \left( U^h_{i-1, j} + U^h_{i, j+1} - 4 U^h_{i, j} + U^h_{i, j-1} + U^h_{i+1, j} \right)                    \\
	=                                                                               & ~ \sum_{m=-1}^1 \sum_{n=-1}^1 U_{i+m,j+n} K_{m,n} \stackrel{\cref{eq:cross-correlation}}{=} \left(U^h \ast K \right)_{ij}
\end{align*}
The equality of the right hand sides follows from~\cref{eq:discrete_f}, showing that a finite difference discretization can be implemented using cross-correlation in a CNN. Next, we will use this analogy to derive a physics-aware loss function in the context of CNNs.

Note that, in this whole subsection, we have neglected the treatment of boundary conditions; we discuss this directly in the context of the application of our approach to the Navier--Stokes equations in~\cref{sec:geometry-bcs}.

\subsection{Derivation of a physics-aware loss function} \label{sec:physics_loss}
Let us consider a generic system of PDEs given in implicit form on the same computational domain $\Omega = [0,1]^2$:
\begin{equation}
	F \left(x,~u\!\left(x\right),~\frac{\partial u}{\partial x_1}\!\left(x\right),~\frac{\partial u}{\partial x_2}\left(x\right), ~\dots  \right) = 0 \text{,} \quad x \in \Omega.
	\label{eq:continuous_pde}
\end{equation}
Here, $F$ is a nonlinear function which may depend on partial derivatives of $u$ of any order. Therefore, \cref{eq:continuous_pde} is a generalization of the diffusion equation~\cref{eq:pde}.

Analogously to~\cref{sec:FD_and_convolutions}, we can discretize~\cref{eq:continuous_pde} by approximating the derivatives using finite differences on a structured $n+1 \times n+1$ grid. 
Appropriate finite differences schemes for various PDEs can be found in the literature; see, for instance,~\cite{Leveque_BasicsFD_2007, smith1985numerical, strikwerda2004finite}. 
In~\cref{sec:ns}, we discuss the specific case of the Navier--Stokes equations, which are the main application discussed in this work. 
As mentioned before, other discretization schemes on structured grids, such as finite element or finite volume discretizations can also be used. 

Let $U^h \in \mathbb{R}^{d \times \left(n+1\right) \times \left(n+1\right)}$
be the tensor representation of the discrete solution with $d$ components. 
Then, the discrete problem corresponding to~\cref{eq:continuous_pde} can be written as
\begin{equation}
	F^h \left(X^h,~U^h\!,~U^h \ast D^x\!,~U^h \ast D^y\!,~\dots\right) = 0,
	\label{eq:discrete_pde}
\end{equation}
where 
$
X^h
=
\left(\mathbf{x}_{ij}^h\right)_{ij}
\in
\mathbb{R}^{2 \times \left(n+1\right) \times \left(n+1\right)}
$
is the tensor containing all the grid nodes, and $D^x$ and $D^y$ are kernel matrices corresponding to the finite difference discretization of the partial derivatives 
$$
	\frac{\partial u}{\partial x_1}
	\quad
	\text{and}
	\quad
	\frac{\partial u}{\partial x_2},
$$
respectively. 
As shown in~\cref{sec:FD_and_convolutions} for the example of a standard five-point stencil, any finite difference discretization of a partial derivative can be written as the cross-correlation with the corresponding finite difference stencil. 
Higher derivatives can therefore be treated analogously.

The generally nonlinear system of equations \cref{eq:discrete_pde} can be reformulated as a least-squares problem for the discrete residual
\begin{equation}
	\arg\min\limits_{U^h} \left\| F^h \left(X^h,~U^h\!,~U^h \ast D^x\!,~U^h \ast D^y\!,~\dots\right) \right\|_2^2.
	\label{eq:min_general}
\end{equation}
Both problems are equivalent when the same boundary conditions are imposed. 
While solving~\cref{eq:discrete_pde} benefits classical numerical solvers, the minimization problem~\cref{eq:min_general} is better suited for a neural network approach. 
The solution tensor $U^h$ on the grid $\Omega_h$ can be replaced by the CNN output $f_{N!N}^\Psi$, resulting in 
\begin{equation}
\arg\min\limits_{\Psi} \left\| F^h \left(X^h,~f_{N\!N}^\Psi\!,~f_{N\!N}^\Psi \ast D^x\!,~f_{N\!N}^\Psi \ast D^y\!,~\dots\right) \right\|_2^2,
\label{eq:pde_loss}
\end{equation}
where $\Psi$ represents the network parameters.
The cross-correlation can be easily implemented since it is a standard operation in state-of-the-art deep learning libraries.

This approach is related to physics-informed neural networks (PINNs) that classically use dense feedforward neural networks for discretization.
In classical PINNs, the residual of the partial differential equation is also minimized in a least-squares sense. 
Therefore, the differential operator is evaluated by automatic differentiation of the network function using the back propagation algorithm~\cite{rumelhart1986learning}. 
Our CNN-based approach differs as it employs a classical discretization, for instance, based on finite difference stencils. 
In our approach, the neural network predicts the coefficients of the discrete solution; cf.~\cref{sec:FD_and_convolutions}.

In \cref{eq:pde_loss}, the CNN input is intentionally omitted.
In fact, if the solution of~\cref{eq:discrete_pde,eq:min_general} is unique due to appropriate boundary conditions and finite difference discretization, the solution does not depend on any input parameters. 
Hence, it is sufficient to train the neural network for a constant output; this could be simply realized via a bias vector in the output layer. 
This approach is not yet relevant in practice, and the discussion in~\cite{GHK:2023:ShortNote} shows that, for a simple stationary diffusion problem, solving~\cref{eq:pde_loss} using an SGD-based optimizer cannot compete with solving the discrete system~\cref{eq:discrete_pde} using classical numerical solvers, such as gradient descent or the conjugate gradient method.

However, the physics-aware loss function~\cref{eq:pde_loss} becomes relevant once the CNN model serves as a surrogate model for multiple configurations parameterized by the input of the CNN model; cf.~\cref{sec:cnn} for the data-based approach. 
Analogously, we consider image representations describing the geometry of the computational domain.

\subsection{Geometry-dependency and boundary conditions} \label{sec:geometry-bcs}

\begin{figure}[t]
	\centering
	\includegraphics[width=\textwidth]{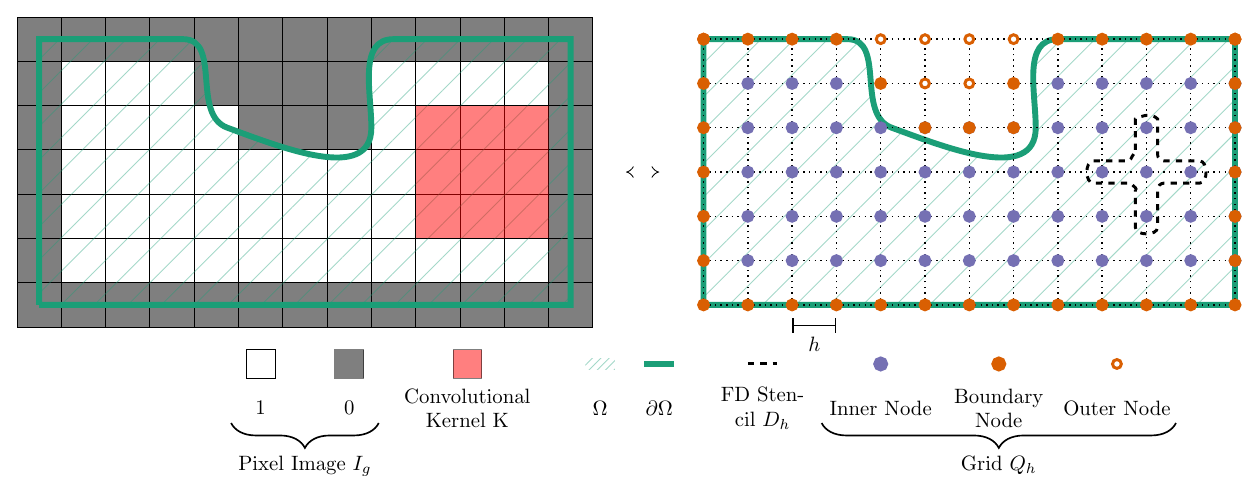}
	\caption{\label{fig:pixel_equiv_grid}A pixel image $I_g$ (left) and the corresponding FD-grid $\Omega_h$ (right) of a geometry $\Omega$, whose border $\partial \Omega$ is drawn in green in the figure.  Applying a five-point finite difference stencil $D_h^k$ to the FD-grid is equivalent to applying a convolutional filter $K$ with fixed weights to the pixel image $I_g$. Note that the values associated with the pixels and their grid node counterparts are not depicted.
	}
\end{figure}

To extend the physics-aware loss function for variations in the geometry of the computational domain, we introduce a discretization of~\cref{eq:continuous_pde} 
that is compatible with the pixel image representation employed in the CNN.
For this purpose, we consider a rectangle $Q$ encompassing $\Omega$ and use an equidistant grid $Q_h$ with grid step size $h$ to discretize it. 
The set of grid nodes is denoted as $X=\{x_{i,j}\}$, with $w$ grid nodes in the $x$ direction and $h$ grid nodes in the $y$ direction; cf. \cref{fig:pixel_equiv_grid}. 

Now, let $I_g$ be an input image describing the geometry of the computational domain $\Omega$; even though we can generally employ any geometry representation, it is important that the boundary pixels are uniquely determined because we explicitly use them in our approach; see, for instance,~\cref{fig:input_images} (left). Plugging the CNN model $f_{N\!N}\!\left( Ig \right)$, as described in~\cref{sec:cnn}, into the physics-aware minimization problem~\cref{eq:pde_loss}, we then obtain
\begin{equation}
	\underset{\Psi}{\arg\!\min} \left\| F^h \left(X^h,~C_{I_g}\left(f_{N\!N}^\Psi \left( I_g \right)\right),~C_{I_g}\left(f_{N\!N}^\Psi \left( I_g \right)\right) \ast D^x\!,~C_{I_g}\left(f_{N\!N}^\Psi \left( I_g \right)\right) \ast D^y\!,~\dots\right) \right\|_2^2.
	\label{eq:single_geometry_pde_loss}
\end{equation}
Here, we only enforce the physics-awareness for those pixels which are inside the computational domain, as indicated by the input image representation $I_g$, and the operator $C_{I_g}$ corresponds to enforcing the boundary conditions.
In the following we will abbreviate
the notation of $F^h$ and write $\left\| F^h \left(X^h, f_{N\!N}^\Psi \left( I_g \right), \dots \right) \right\|$ for ease of readability.

In order for~\cref{eq:single_geometry_pde_loss} to be well-defined, we have to prescribe boundary conditions at the boundary pixels. There are at least two ways of enforcing boundary conditions in the context of physics-based neural network models. In particular, we can either add a loss term associated with the boundary conditions or explicitly encode the boundary conditions in the network function. In the literature, the former is is also denoted as \textit{soft enforcement} of boundary conditions, whereas the latter is denoted as \textit{hard enforcement} of boundary conditions; cf. \cite{SGPW:2020:SMF} for a more detailed discussion. It has been observed in the literature that soft enforcement of boundary conditions can be problematic in different ways: it can make the training less robust and also may lead to cases where the training does not converge to the solution; see for example \cite{SGPW:2020:SMF}. Therefore, we focus on hard enforcement of boundary conditions. For instance, in case of Dirichlet boundary conditions, can be easily done by explicitly writing the correct values in the output image of the neural network before applying the loss function; at the same time, and as in classical discretization methods, we do not enforce the physical loss in those pixels. 
The explicit enforcement of the boundary conditions is indicated by the operator $C_{I_g}$ in~\cref{eq:single_geometry_pde_loss}. If different boundary conditions are prescribed on different parts of the boundary, we encode this by specific values in the input image; cf.~\cref{fig:input_images}. See also~\cref{subsec:boundary_conditions} for a specific discussion of our implementation of boundary conditions for the Navier--Stokes equations.

Now, we extend the training of the surrogate model to multiple geometries. Therefore, analogously to the data-driven case~\cref{eq:mse}, we optimize the loss function over a training data set of geometries $T$, resulting in the following loss function:
\begin{equation}
		\underset{f_{N\!N}}{\arg\!\min} \frac{1}{\vert T \vert} \sum_{g \in T}	\left\| F^h \left(X^h, f_{N\!N}^\Psi \left( I_g \right), \dots \right) \right\|_2^2.
	\label{eq:multiple_geometry_pde_loss}
\end{equation}
The main difference to the data-based loss function~\cref{eq:mse} from~\cref{sec:cnn} is that no reference flow data $f_g$ but only the mathematical model of the PDE is necessary for the training. Of course, both loss functions can also be combined into a hybrid loss function
\begin{equation}
		\underset{\Psi}{\arg\!\min} \frac{1}{\vert T \vert} \sum_{g \in T}
		\omega_{\text{PDE},g} \left\| F^h \left(X^h, f_{N\!N}^\Psi \left( I_g \right), \dots \right) \right\|_2^2 + \omega_{\text{data},g} \left\| f_{N\!N}^\Psi \left( I_g \right) - u_g \right\|^2,
	\label{eq:multiple_geometry_hybrid_loss}
\end{equation}
where the weights $\omega_{\text{PDE},g}$ and $\omega_{\text{data},g}$ balance the two loss terms. Different variants of the hybrid loss function are possible, for instance, 
$$
\begin{array}{lcll}
	\omega_{\text{data},g} = \alpha & \land & \omega_{\text{PDE},g} = 0, & \text{if reference data is available}, \\
	\omega_{\text{data},g} = 0 & \land & \omega_{\text{PDE},g} = \beta, & \text{otherwise}, \\
\end{array}
$$
or
$$
\begin{array}{lcll}
	\omega_{\text{data},g} = \alpha & \land & \omega_{\text{PDE},g} = \beta, & \text{if reference data is available}, \\
	\omega_{\text{data},g} = 0 & \land & \omega_{\text{PDE},g} = \beta, & \text{otherwise}. \\
\end{array}
$$
Here, $\alpha, \beta > 0$ are some weight parameters. Other strategies for choosing the weights are, of course, also possible. For a theoretical discussion based on the neural tangent kernel on how to balance PDE and data loss terms for classical PINNs, see~\cite{wang_when_2022}. 

Next, we discuss the details of our model for the specific problem considered here, that is, the Navier--Stokes equations~\cref{eq:ns}.

\section{Application to the Navier--Stokes equations} \label{sec:ns}

In this work, we are concerned with the application of our approach to the Navier--Stokes equations. Therefore, we discuss, in this section, the derivation of the physics-aware loss function and the treatment of the boundary conditions.

\subsection{Physics-aware loss function}
We have already introduced the Navier--Stokes equations in~\cref{eq:ns} of~\cref{sec:flow-problem}. If we expand it in terms of the individual components, we obtain
\begin{align}
	u \frac{\partial u}{\partial x} + v \frac{\partial u}{\partial y} + \frac{\partial p}{\partial x} - \nu \left( \frac{\partial^2 u}{\partial x^2} + \frac{\partial^2 u}{\partial y^2} \right) &= 0 \label{eq:component_moment_1}\\
	u \frac{\partial v}{\partial x} + v \frac{\partial v}{\partial y} + \frac{\partial p}{\partial y} - \nu \left( \frac{\partial^2 v}{\partial x^2} + \frac{\partial^2 v}{\partial y^2} \right) &= 0 \label{eq:component_moment_2} \\
	\frac{\partial u}{\partial x} + \frac{\partial v}{\partial y} &= 0 \label{eq:component_div} \text{,}
\end{align}
where $u$ and $v$ are the $x$- and $y$-component of the flow field. 
Here, \cref{eq:component_div} is the continuity equation and \cref{eq:component_moment_1} and \cref{eq:component_moment_2} are the components of the momentum equation. 

We discretize~\cref{eq:component_moment_1,eq:component_moment_2,eq:component_div} using the central-difference stencils 
$$
\begin{aligned}
	\left(D^{x} u_h\right)_{i+j}  := \frac{u_h^{i+1,j} - u_h^{i-1,j}}{2h}, \quad 
	&
	\left(D^{xx} u_h \right)_{i+j} := \frac{u_h^{i+1,j} - 2u_h^{i,j} + u_h^{i-1,j}}{h^2}, \\
	\left(D^{y} u_h\right)_{i+j}  := \frac{u_h^{i,j+1} - u_h^{i,j-1}}{2h}, \quad
	&
	\left(D^{yy} u_h \right)_{i+j} := \frac{u_h^{i,j+1} - 2u_h^{i,j} + u_h^{i,j-1}}{h^2}. \\
\end{aligned}
$$ 
 The resulting discretized Navier--Stokes equations read
\begin{align}
	u_h \odot D^{x} u_h + v_h \odot D^{y} u_h  + D^{x} p_h - \nu \left(D^{xx} u_h + D^{yy} u_h \right) &= 0 \label{eq:discrete_component_moment_1}\\
	u_h \odot D^{x} v_h + v_h \odot D^{y} v_h  + D^{y} p_h - \nu \left(D^{xx} v_h + D^{yy} v_h \right) &= 0 \label{eq:discrete_component_moment_2} \\
	D^{x} u_h + D^{y} v_h &= 0, \label{eq:discrete_component_div}
\end{align}
where $\odot$ is the Hadamard product, that is, element-wise product. Following the discussion in~\cref{sec:FD_and_convolutions}, this can equivalently be written using the cross-correlation  operation $\ast$ as follows:
\begin{align}
	U_h \odot U_h \ast D^{x} + V_h \odot U_h \ast D^{y} + P_h \ast D^{x} - \nu \left( U_h \ast D^{xx} + U_h \ast D^{yy} \right) & = 0 \label{eq:discrete_corr_component_moment_1}\\
	U_h \odot V_h \ast D^{x} + V_h \odot V_h \ast D^{y} +  P_h \ast D^{y} - \nu \left( V_h \ast D^{xx} +  V_h \ast D^{yy} \right) &= 0 \label{eq:discrete_corr_component_moment_2} \\
	U_h \ast D^{x} + V_h \ast D^{y} &= 0, \label{eq:discrete_corr_component_div}
\end{align}
where, for simplicity, we overload the notation for the discrete differential operators $D^{x}$, $D^{xx}$, $D^{y}$, and $D^{yy}$ with the corresponding stencil matrices, and $U_h$, $V_h$, and $P_h$ are the matrix representations of the solution fields corresponding to $u_h$, $v_h$, and $p_h$, respectively. Furthermore, for simplicity, we omit the treatment of the boundary conditions for now and refer to~\cref{subsec:boundary_conditions} for a detailed discussion.

\Cref{eq:discrete_corr_component_moment_1,eq:discrete_corr_component_moment_2,eq:discrete_corr_component_div} correspond to the discrete nonlinear system of equations
\begin{align}
	& N(U_h,V_h) + G(P_h) & = 0 \quad & & \text{in } \Omega, \label{eq:operator_moment} \\
	& D(U_h,V_h) & = 0 \quad & & \text{in } \Omega, \label{eq:operator_div}
\end{align}
where each of the operators $N$, $G$, and $D$ can be implemented using the cross-correlation and the Hadamard product as the building blocks; cf.~\cref{eq:discrete_corr_component_moment_1,eq:discrete_corr_component_moment_2,eq:discrete_corr_component_div}. The operator $N$ is nonlinear, whereas $G$ and $D$ are both linear operators.

Now, we apply the physics-aware surrogate modeling approach described in~\cref{sec:geometry-bcs} to predict the solution of~\cref{eq:operator_moment,eq:operator_div} for varying geometries. The surrogate model takes the form of
\begin{align*}
	f_{N\!N}^\Psi: \mathbb{R}^{w \times h} &\rightarrow \mathbb{R}^{3\times w\times h} \\
	I_g &\rightarrow \begin{pmatrix}
		U_{N\!N}^\Psi(I_g) \\
		V_{N\!N}^\Psi(I_g) \\
		P_{N\!N}^\Psi(I_g)
	\end{pmatrix}\text{,}
\end{align*}
where $I_g$ is, again, the pixel image representation of a geometry $g$, and $U_{N\!N}^\Psi\left(I_g\right)$, $V_{N\!N}^\Psi\left(I_g\right)$, and $P_{N\!N}^\Psi\left(I_g\right)$ correspond to the CNN predictions for the matrices resp.~images $U_h$, $V_h$, and $P_h$, respectively.

Combining our physics-aware approach for multiple geometries as described in~\cref{sec:geometry-bcs} with this CNN model and the discrete residual of the Navier--Stokes equations~\cref{eq:operator_moment,eq:operator_div}, we obtain the loss function
\begin{equation}
	\frac{1}{\vert T \vert} \sum_{g \in T} \left( \omega_{M}\Vert N(U_{N\!N}^\Psi(I_g),V_{N\!N}^\Psi(I_g)) + G(P_{N\!N}^\Psi(I_g)) \Vert_2^2 + \omega_{D} \Vert D(U_{N\!N}^\Psi(I_g),V_{N\!N}^\Psi(I_g))\Vert_2^2 \right).
	\label{eq:ns_pde_loss}
\end{equation}
Here, $\omega_{M}$ and $\omega_{D}$ are weights for the two loss terms, and $T$ is, again, the set of all training geometries $g$.

To complete our discussion of the application of the physics--aware approach to the Navier--Stokes equations, we discuss the specific treatment of the boundary conditions in the next section.

\subsection{Treatment of boundary conditions}
\label{subsec:boundary_conditions}

\begin{figure}[t]
	\centering
	\begin{subfigure}[c]{0.49\textwidth}
		\includegraphics[width=\textwidth]{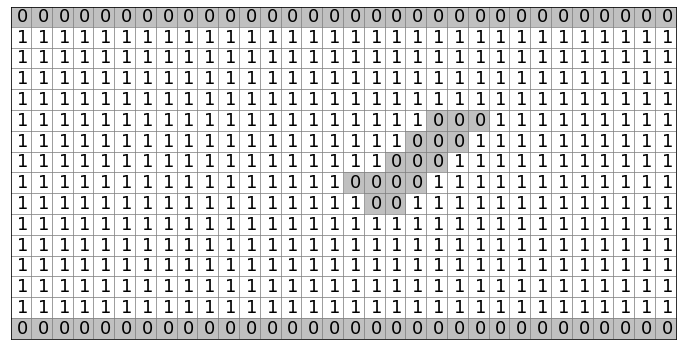} 
		\subcaption{\label{fig:geometry_input_image}Geometry image}
	\end{subfigure}
	\begin{subfigure}[c]{0.49\textwidth}
		\includegraphics[width=\textwidth]{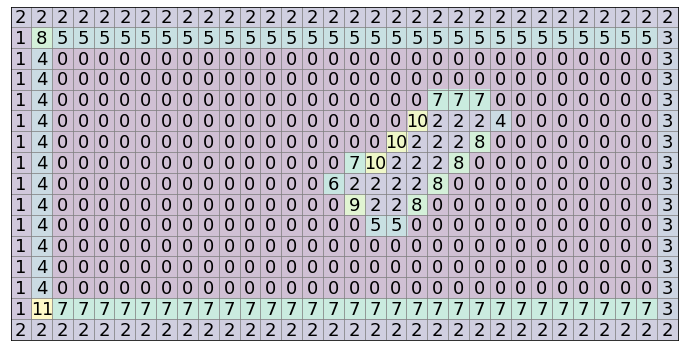} 
		\subcaption{\label{fig:boundary_input_image}Boundary image}
	\end{subfigure}
	\caption{Low-resolution pixel image inputs that are used by our model. The geometry image \cref{fig:geometry_input_image} is passed as input to the CNN and the boundary image \cref{fig:boundary_input_image} is used for the construction of the physics-aware loss.  The geometry represented here was previously described in \cref{fig:showcase_domain}
		\label{fig:input_images}
	}
\end{figure} 

As discussed in~\cref{sec:geometry-bcs}, we enforce Dirichlet boundary conditions explicitly by hard-coding the values of the pixels in the output image. In particular, for our boundary value problems, as introduced in~\cref{sec:flow-problem}, we consider the following boundary conditions: For inlet boundary condition, we set 
\begin{equation} \label{eq:hard_inflow_constraint}
	U_{1,j} = 3, V_{1,j} = 0, \quad \forall j = 1,\ldots,h-1,
\end{equation}
where $h$ (height) is the number of pixels in $y$ direction. Moreover, the no-slip boundary conditions
\begin{equation} \label{eq:hard_no_slip_constraint}
	U_{i,1} = U_{i,h} = V_{i,1} = V_{i,h} = 0, \quad \forall i = 1,\ldots,w,
\end{equation}
are enforced at the lower and upper walls as well as the zero pressure boundary condition 
\begin{equation} \label{eq:hard_outflow_constraint} 
	P_{w,j} = 0, \quad \forall j = 1,\ldots,h-1,
\end{equation}
at the outlet. Here, $w$ (width) is the number of pixels in $x$ direction.

Neumann pressure boundary conditions can be implemented by introducing ghost-nodes outside our computational domain.
However, in this work we do not consider Neumann boundary conditions and instead refer to \cite{GSW:2021:PGN}.
Note, though, that the use of ghost-nodes may not be feasible in the case of irregular obstacle boundaries.
Here, interpolation to a pixel image may introduce corner points on the boundary where the normal vector is not well defined.
	
Further, to avoid the usage of pressure values in pixels where the pressure is not defined, we employ one-sided differences in pixels adjacent to corresponding boundaries.
We encode the different stencils in an additional input image; cf.~\cref{fig:input_images} (right). 
The numbering scheme used for the grid nodes is as follows: $0$ represents internal nodes, $1$ corresponds to inflow boundary nodes, $2$ represents no-slip boundary nodes, and $3$ denotes outflow boundary nodes. Nodes with numbers $4$ and above require one-sided approximations for the pressure gradient.
It is important to note that some pixels correspond to nodes outside the original domain $\Omega$ due to obstacles. These nodes are marked as $0$ in the geometry image and $2$ in the boundary image. Velocity and residual values are set to $0$ in these nodes. This is necessary as the governing equations are not defined in those nodes.

\begin{figure}[t]
	\centering
	\begin{subfigure}[c]{\textwidth}
		\includegraphics[width=\textwidth]{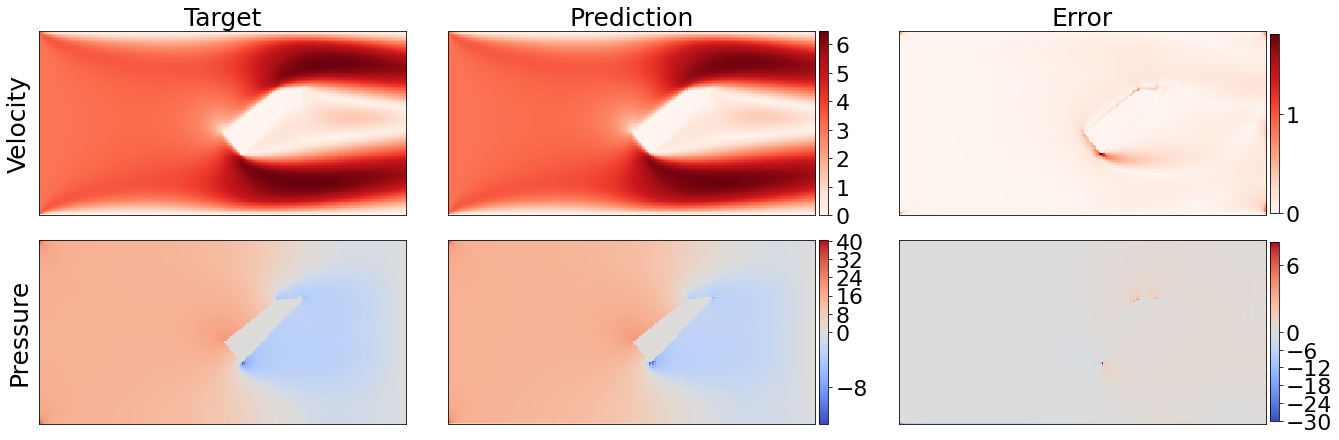} 
		\caption{\label{fig:single_geometry_1896}Velocity and pressure for the physics-aware approach (Prediction) compared to the OpenFOAM simulation on the \textit{locally refined} mesh (Target).}
	\end{subfigure}
	
	\begin{subfigure}[c]{\textwidth}
		\includegraphics[width=\textwidth]{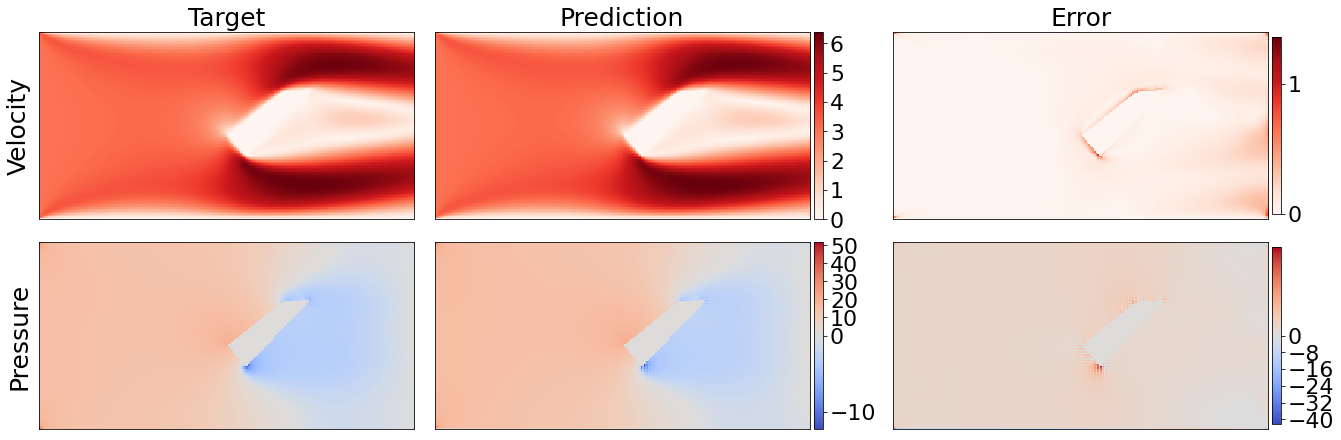} 
		\caption{\label{fig:single_geometry_1896_rasterized}Velocity and pressure for the physics-aware approach (Prediction) compared to the OpenFOAM simulation on the \textit{Cartesian rasterized} mesh (Target).}
	\end{subfigure}
	\caption{\label{fig:single_geometry} Results for a single geometry. CNN model with Swish activation function and $5 \cdot 10^{-5}$ as the learning rate for the Adam optimizer.}
\end{figure}

\subsection{Example on a single geometry} \label{subsec:single-geometry}
In order to verify that the physics-aware loss enables us to learn a solution of the Navier--Stokes equations, we consider a single fixed geometry. This is not useful in practice since we could more efficiently directly discretize the Navier--Stokes equations using finite differences and solve the discrete system using suitable numerical solvers; see~\cite{GHK:2023:ShortNote} for a comparison for a simple Laplace problem. 

We compare the model prediction against FVM simulations with OpenFOAM on two different meshes: a locally refined mesh (\cref{fig:showcase_mesh}) and using the same pixel grid as the CNN model (removing the pixels inside the obstacle). The results are plotted in~\cref{fig:single_geometry}. Compared with the simulation on a locally refined mesh, we obtain low relative $L_2$ errors (defined in~\cref{eq:rel_l2_error}) of $2.6\,\%$ for $u$ and $2.8\,\%$ for $p$. As can be seen in~\cref{fig:single_geometry_1896}, the velocity error is particularly high near the obstacle, presumably due to non-resolved boundary layers. The comparison against the simulation on the rasterized mesh in~\cref{fig:single_geometry_1896_rasterized} shows a visual improvement of these errors, and the relative $L_2$ error for the velocity reduces to $2.2\%$. This suggests that part of the error is due to insufficient mesh resolution. An error of $0$ cannot be obtained since the CNN model is based on a finite difference discretization whereas the reference data is computed using FVM simulations for both types of meshes.

In total, we conclude from the results for a single geometry that a CNN model with physics-aware loss may learn a good approximation of the solution of the Navier--Stokes equations. Later, in~\cref{sec:results}, we will investigate the performance of the CNN-based surrogate model trained on a data set consisting of multiple geometries, introducing another level of complexity to the model.	

\section{Architecture of the convolutional neural network} \label{sec:model}
In this section, we discuss the network architecture of our CNN-based surrogate models, utilizing the same architecture type for both the data-based and physics-aware models described in~\cref{sec:cnn,sec:method} respectively.
We employ a fully convolutional neural network that only performs convolutions, up- or downsampling.
For a comprehensive understanding of CNNs, we refer to~\cite[Chapt. 9]{goodfellow2016deep} and the references therein.

Our CNN architecture draws inspiration from the U-Net architecture~\cite{RFB:2015:UNET}. 
It consists of an encoder, transforming input image(s) into a lower dimensional representation in the bottleneck, and a decoder, transforming the bottleneck output into velocity and pressure output images. 
The U-Net architecture's symmetric encoder and decoder paths are connected via skip connections.
The performance of data-based surrogate models with U-Net architecture is generally superior compared to bottleneck CNNs without skip connections, as discussed in~\cite{eichinger_surrogate_2022}.

\begin{figure}[t]
	\centering
	\begin{subfigure}[c]{\textwidth}
		\includegraphics[width=\textwidth]{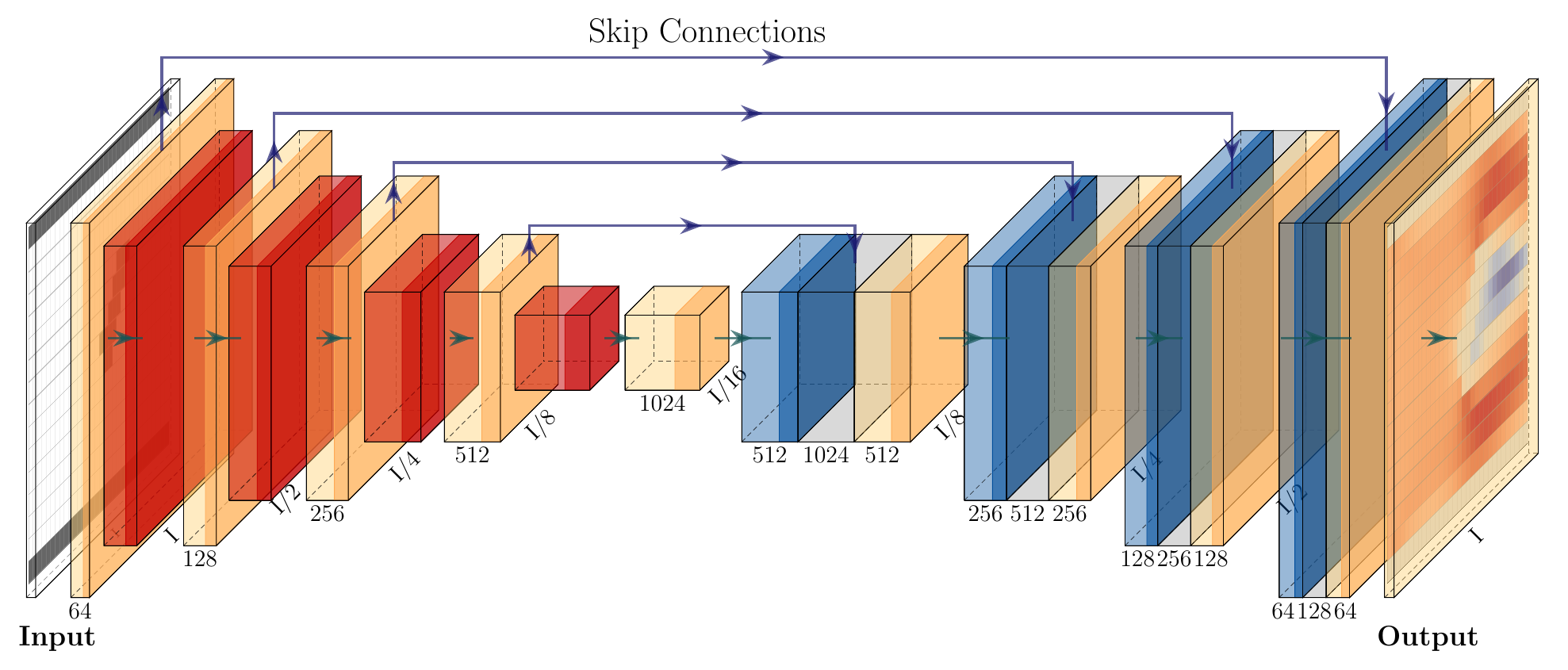} 
		\subcaption{Model}
	\end{subfigure}
	\begin{subfigure}[c]{\textwidth}
		\includegraphics[width=\textwidth]{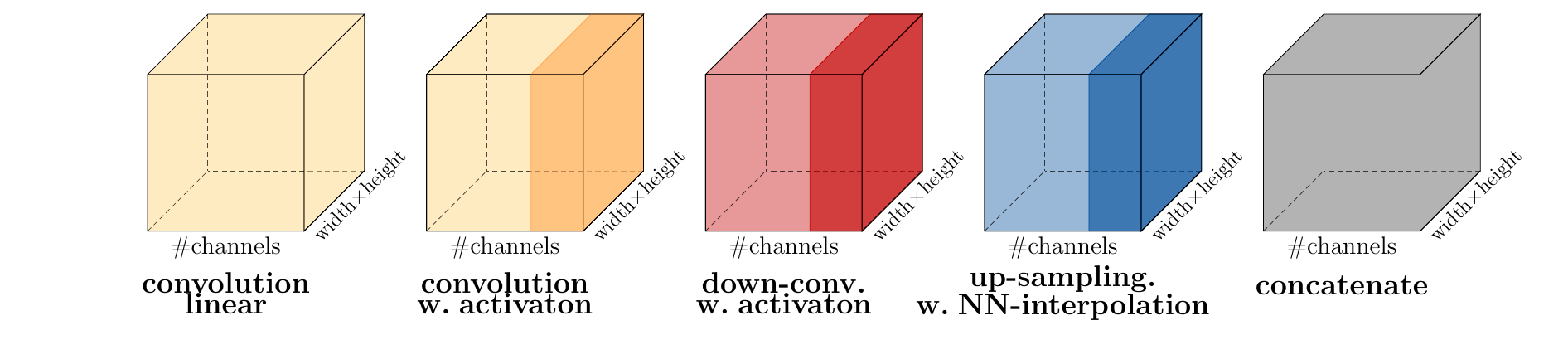} 
		\subcaption{Legend}
	\end{subfigure}
	\caption{\label{fig:unet} Exemplary model architecture with four levels.}
\end{figure} 

The encoder of our network consists of blocks comprising a $3 \times 3$ convolutional layer, followed by a $2 \times 2$ convolutional layer with a stride of 2, both followed by an activation. 
The first convolution extracts input features, while the second convolution reduces spatial dimensions. 
We avoid max pooling for downsizing to prevent high-frequency artifacts as discussed in~\cite{HS:2016:GLP}.
The decoder mirrors the encoder and includes upsampling layers with nearest-neighbor interpolation followed by a convolutional layer with an activation. 
This upsampling technique helps avoiding checkerboard artifacts mentioned in~\cite{ODO:2016:DCA} that can occur with deconvolutional and downward convolutional layers. 
Matching encoder and decoder blocks are connected by skip connections following the U-Net architecture. 
The encoder block output is concatenated with the upsampling layer output, doubling the number of filters. 
Refer to~\cref{fig:unet} for an illustration of this architecture with one decoder path and four levels. In our experiments in~\cref{sec:results}, unless stated otherwise, we employ an 8-level model. 
Additionally, it should be noted that our models use separate decoder paths for each scalar output field, $U$, $V$, and $P$.

We refer to~\cref{sec:hyperparams} for additional comments on the choice of hyper parameters, including the model architecture.

\section{Generation of (training) data} \label{sec:trainingdata}

\begin{figure}[t]
	\centering
	\begin{subfigure}[c]{0.56\textwidth}
		\includegraphics[width=\textwidth]{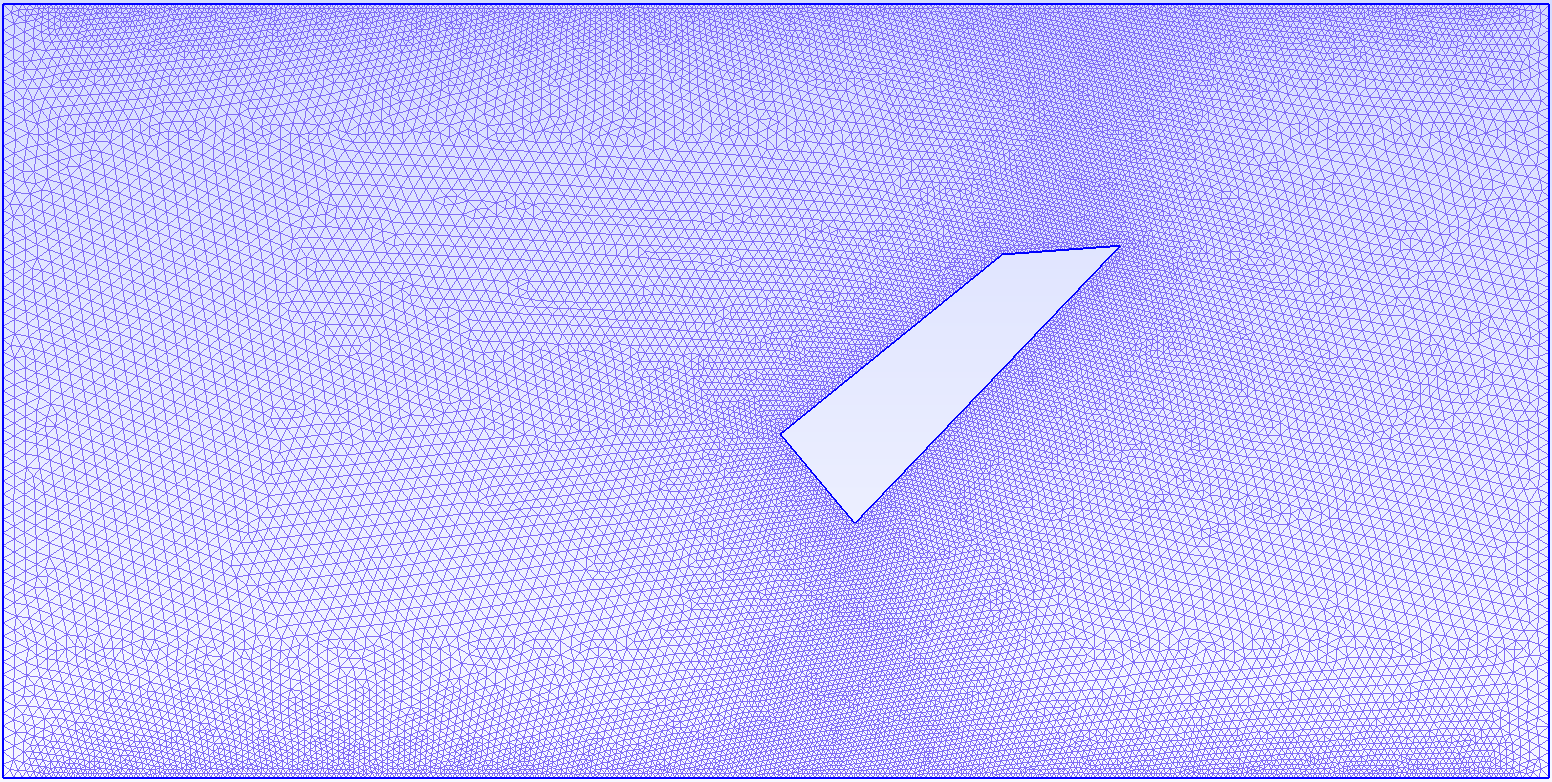} 
		\caption{\label{fig:showcase_mesh}}
	\end{subfigure}
	\begin{subfigure}[c]{0.4\textwidth}
		\includegraphics[width=\textwidth]{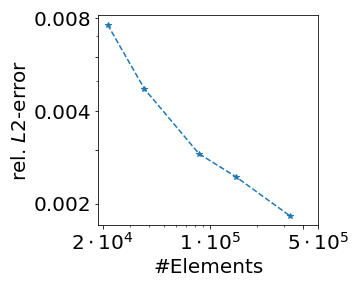} 
		\caption{\label{fig:mesh_convergence}}
	\end{subfigure}
	\caption{(a) Example of a locally refined mesh used for simulations. (b) Mesh Convergence plot with relative errors compared against on a reference simulation for a fine mesh with $\approx 1\,300\,000$ elements.
		The depicted mesh corresponds to the second node in the convergence plot.
	}
\end{figure}

\begin{figure}[t]
	\centering
	\includegraphics[width=0.8\textwidth]{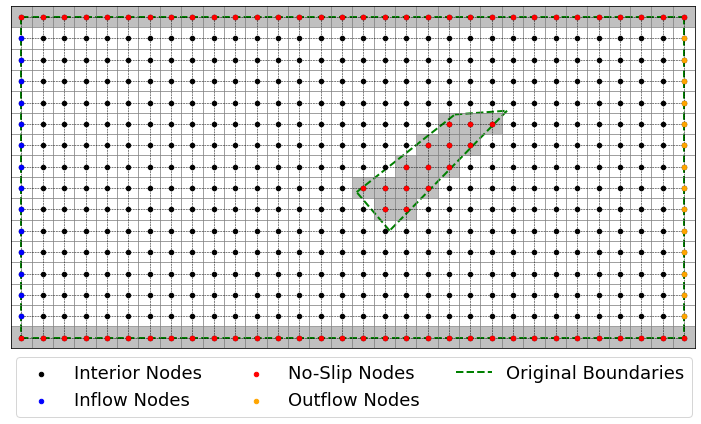} 
	\caption{\label{fig:showcase_interpolation} Exemplary representations of the interpolation process, here with a lower resolution of $32\times 16$. The green dotted line shows the border of the computational domain $\Omega_P$. Note that the original boundary is plotted behind the outer nodes.}
\end{figure}

\begin{figure}[t]
	\centering
	\begin{subfigure}[c]{0.49\textwidth}
		\includegraphics[width=\textwidth]{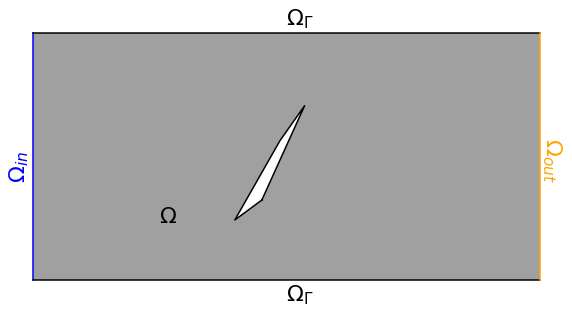} 
		\caption{\label{fig:showcase_sharp_angle_geometry}}
	\end{subfigure}
	\begin{subfigure}[c]{0.49\textwidth}
		\includegraphics[width=\textwidth]{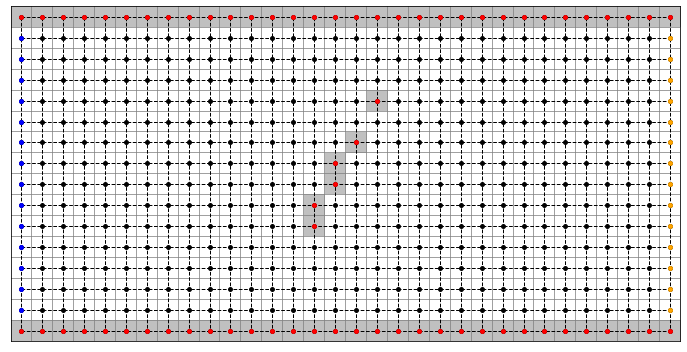} 
		\caption{\label{fig:showcase_sharp_angle_interpolation}}
	\end{subfigure}
	\caption{\label{fig:showcase_sharp_angle_geo_and_int} Exemplary representations of gross distortion caused by too acute angles. (a) The original geometry and (b) the reduced pixel image, here with a lower resolution of $32\times 16$. The green dotted line shows the borders of the computational domain $\Omega_P$. Note that there the original boundary is plotted behind the outer nodes.
	}
\end{figure}

\begin{figure}[t]
	\centering
	\begin{subfigure}[c]{0.49\textwidth}
		\includegraphics[width=\textwidth]{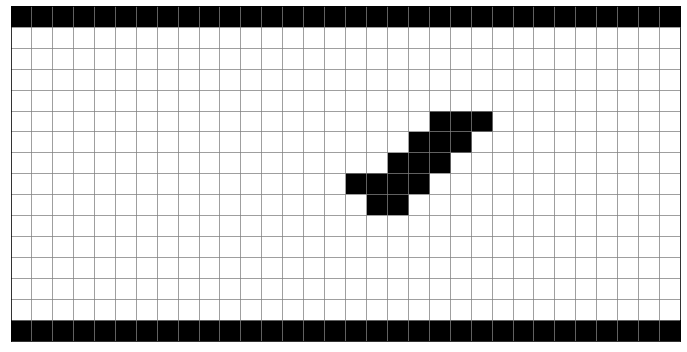} 
		\subcaption{Geometry}
	\end{subfigure}
	\begin{subfigure}[c]{0.49\textwidth}
		\includegraphics[width=\textwidth]{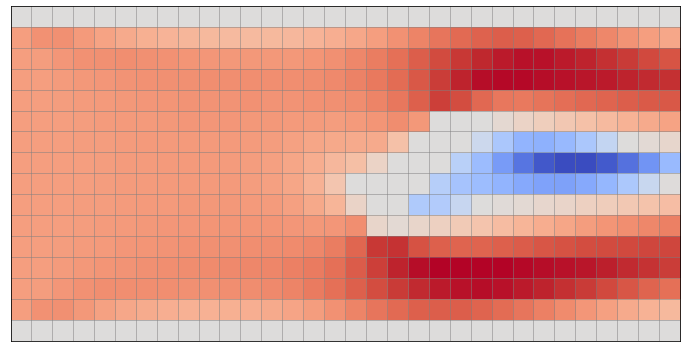} 
		\subcaption{$x$-Velocity}
	\end{subfigure}
	
	\begin{subfigure}[c]{0.49\textwidth}
		\includegraphics[width=\textwidth]{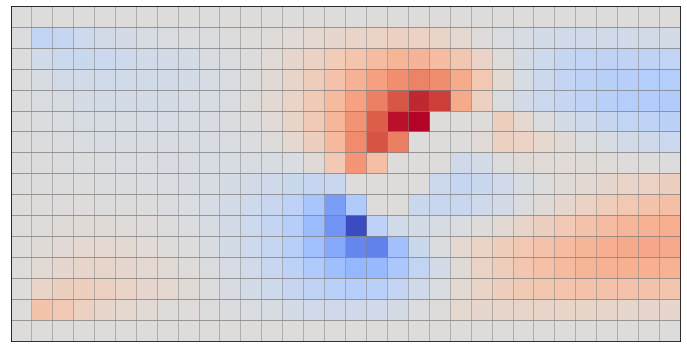} 
		\subcaption{$y$-Velocity}
	\end{subfigure}
	\begin{subfigure}[c]{0.49\textwidth}
		\includegraphics[width=\textwidth]{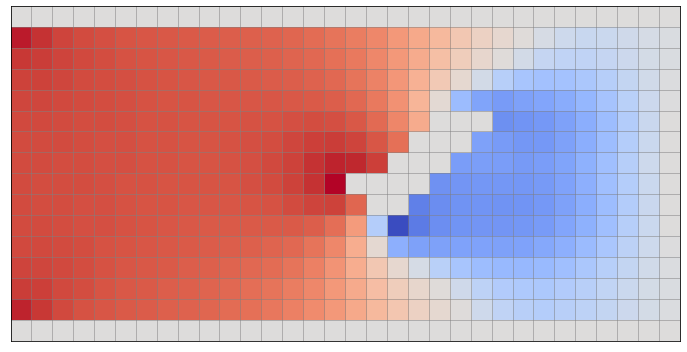} 
		\subcaption{Pressure}
	\end{subfigure}
	\caption{\label{fig:showcase_data} (a)-(d) Exemplary representations of the pixel images, here with a lower resolution of $32\times 16$.}
\end{figure}

Training our surrogate models for the challenging task of predicting solutions across a wide range of geometries necessitates a substantial dataset. 
For the fully data-driven model, reference data is required along with input images of the geometries. 
In contrast, the physics-aware model solely relies on the input images, as the physics loss replaces the need for reference data. 
Thus, generating training data involves creating random obstacle geometries and their pixel image representations. 
To validate our model, we generate simulation meshes and conduct corresponding CFD simulations as reference data.

\subsection{Input data} \label{sec:input}
The input data is needed for the training of the data-driven and physics-aware model as well as of hybrid variants.

\paragraph{Geometry generation} 
In line with~\cref{sec:flow-problem}, our focus lies on two-dimensional rectangular channel geometries featuring star-shaped obstacles that do not touch the boundary (inlet, outlet, upper and lower walls). 
The computational domain is given by $\Omega = \left[ 0, 6 \right] \times \left[ 0, 3 \right]$ from which we have removed a star-shaped obstacle $P$ defined by its corners. 
The obstacle's corners are randomly positioned around a central point, inspired by the approach outlined in~\cite{EHK:2021:FPN}. 
We only consider obstacles with a maximum width of $50\%$ of the channel's height and a minimum distance of $0.75$ from any boundary. 
To prevent significant distortion of the geometry in the image representation, we impose a minimum angle of 10° at each vertex of the obstacle. 
This prevents excessively acute angles, as shown in~\cref{fig:showcase_sharp_angle_geometry}, which could result in a disconnected obstacle representation.

\paragraph{Geometry image representation} 
As discussed before, CNNs rely on input data with a tensor-product structure. 
Therefore, we interpolate the geometry to a binary $256 \times 128$ pixel image. 
In~\cref{fig:showcase_interpolation}, an exemplary geometry is shown in its pixel image representation, where white pixels (encoded as $1$) correspond to fluid cells, while gray pixels (encoded as $0$) represent walls and the obstacle. The pixel value is determined based on whether the center of the pixel is within the fluid domain, or not.
In previous studies, signed distance function (SDF) representations were also used to describe the geometry, generally leading to slightly better results; cf.~\cite{EHK:2021:FPN,eichinger_surrogate_2022, GLI:2016:NFA}.
For simplicity, we restrict ourselves to binary input images, as they are close to practically relevant cases; for instance, binary images can be directly generated from imaging techniques, such as magnetic resonance imaging (MRI).

\subsection{Reference output data} \label{sec:reference_ouput_data}

Reference output data is required in order to compute the data loss~\cref{eq:mse}, which is required for the data-driven and hybrid modeling approaches. In case of the fully physics-aware approach, the reference data is only used for validation. 

\paragraph{Mesh generation} For each obstacle chosen as discussed in~\cref{sec:input}, we generate a computational mesh using Gmsh~\cite{GR:2009:GMSH}; we refer to its documentation\footnote{\href{https://gmsh.info/doc/texinfo/gmsh.html}{https://gmsh.info/doc/texinfo/gmsh.html}} for details on the mesh generation. 
In particular, we use the Frontal Delaunay algorithm~\cite{rebay1993efficient} to generate an unstructured triangular mesh for each case. 
Then, we refine the mesh near all walls, that is, near the upper and lower wall as well as near the obstacle, to resolve the boundary layers of the flow. 
In order to find a suitable level of refinement for the simulations, we performed a mesh convergence study on a representative geometry; cf.~\cref{fig:mesh_convergence}. 
Note that we created a new mesh for each level of refinement, rather than refining an existing mesh.
\Cref{fig:showcase_mesh} shows an mesh with refinement near the boundaries with $\approx 40\,000$ elements, whereas we used meshes with $\approx 160\,000$ -- $200\,000$ to generate the our simulations; we do not display such a mesh for the sake of clarity.

\paragraph{CFD simulations}
The CFD simulations for the generation of the reference data have been performed using OpenFOAM v8~\cite{G:2017:OF}, a software based on the finite volume method (FVM). 
We utilized the simpleFoam solver, which solves the stationary incompressible Navier--Stokes equations
using the semi-implicit method for pressure linked equations (SIMPLE) algorithm \cite{patankar1983calculation,FPS:2020:CMFD}. 
We refer to the OpenFOAM documentation\footnote{\href{https://doc.cfd.direct/openfoam/user-guide-v8/}{https://doc.cfd.direct/openfoam/user-guide-v8/}} for more details on the simpleFoam solver.
The configuration is based on the pitzDaily example, adapted to laminar flow and with stricter convergence criteria.

\paragraph{Interpolation of the simulation data} 
To compare the surrogate model's output with the reference simulation data, we interpolate the simulation data onto the same pixel grid. 
This involves evaluating the FVM solution at the centroids of the pixels, resulting in pixel images $U_h$, $V_h$, and $P_h$ representing the velocity in the $x$ and $y$ directions and the pressure, respectively.

It is important to note that values outside the computational domain $\Omega_p$ are explicitly set to $0$ in both the simulation and the model prediction. This is clearly visible in our plots, such as in~\cref{fig:showcase_data}. To achieve this, we mask the output images based on the geometry representation in the input image.

\section{Some comments on hyperparameter choices} \label{sec:hyperparams}
Our surrogate models depend on numerous hyperparameters, and and their specific choices can significantly affect performance.
These hyperparameters encompass model architecture, such as the number of channels per convolutional layer, the depth of the U-net architecture, and the activation function.
They also include optimizer parameters like the learning rate, learning rate schedule, and batch size. 
Additionally, there are hyperparameters related to the loss function, such as the weights assigned to individual loss terms and discretization parameters for the physics-aware loss, including the employed FD stencils.
Furthermore, there are hyperparameters in a broader sense, for instance, the resolution of input and output images.

Considering the large number of hyperparameters, an exhaustive investigation of their impact is impractical. 
Therefore, instead of conducting a comprehensive grid-search, we have fixed some hyperparameters while varying individual ones. 
To maintain brevity, we provide qualitative discussion of the outcomes rather than presenting extensive results for this process.

\paragraph{Network architecture}
Hyperparameters related to the model architecture determine the number of parameters $\Psi$ and hence the model's capacity to approximate the solution operator; cf.~\cref{sec:cnn}. Due to the high complexity of the solution operator, we expect that the model requires a large number of parameters, whereas a too large number of parameters may lead to overfitting. 
We individually optimize hyperparameters for the network architecture with regards to the validation errors; note that the errors are computed with respect to the reference simulation data.
For the configurations described in~\cref{sec:input} and tested in~\cref{sec:results}, we have determined that a depth of $8$ levels and $64$ channels in the first layer yield a good compromise:
lower values result in reduced approximation properties of the model, while higher values lead to an increased computational effort and decreased generalization properties (in terms of the validation error).

\paragraph{Activation function and learning rate}
As the activation function we either use the rectified linear unit (ReLU)~\cite{PMLR:2011:RELU} or the Swish function~\cite{RZL:2017:SAF}. Whereas for the single geometry case discussed in~\cref{subsec:single-geometry} the use of the Swish activation function was beneficial, the ReLU function generally led to better results for training a surrogate model for multiple geometries; cf.~\cref{sec:results}.
In combination with ReLU we always achieved the best results with a learning rate of $1e-4$.
With swish the optimal learning rate depended on the considered geometry and ranged from $1e-5$ to $1e-4$.

\paragraph{Optimizer and batch size}
The best results for our surrogate model were obtained by using the stochastic gradient descent optimizer with adaptive moment estimation (Adam)~\cite{Kingma:2014:MSO} and a batch size of $1$.
Other optimizers were unable to reliably find suitable minima and greater batch sizes led to greatly increased errors.

\paragraph{Image resolution}
Our fully CNN model architecture, can be applied to any resolution with a power of $2$ number of pixels in both the $x$ and $y$ directions.
The number of pixels in each direction does not have a systematic dependence.
However, if the resolution is too low, the geometry representation may be inaccurate, and the FD discretization error could be high. 

Conversely, higher image resolutions may lead to increased computational effort and reduced accuracy due to limited model capacity. 
We have observed that models for higher resolutions require increased depth to achieve meaningful predictions. 
In addition, training larger and deeper models presents a more difficult optimization problem. We have found that models for pixel images with a width of $512$ pixels and a height of $256$ pixels and larger do not converge to suitable minima as reliably as models for smaller pixel images.
Based on these considerations, we have chosen an image resolution of $256$ pixels in width and $128$ pixels in height, which has also been used in previous studies; cf. \cite{GLI:2016:NFA,EHK:2021:FPN}.

\section{Computational results} \label{sec:results}
In this section, we present numerical results for our surrogate modeling approach.
We investigate the fully data-driven approach (\cref{subsec:supervised_approach}), the fully physics-aware approach (\cref{subsec:unsupervised_approach}), and the hybrid approach (\cref{subsec:mixed_approach}), which combines both.

We evaluate the models using on the relative $L_2$-errors
\begin{equation} \label{eq:rel_l2_error}
	\frac{\Vert U_{N\!N} - U_h \Vert_2}{\Vert U_h \Vert_2}
\end{equation}
as the performance measure. Here, $U_{N\!N}$ is the prediction of our model, and $U_h$ is the reference solution; unless otherwise stated, the reference solution corresponds to the result of an OpenFOAM simulation on a locally refined mesh evaluated at the midpoints of the pixels; cf.\ the discussion in~\cref{sec:trainingdata}. We then compute the $L_2$-error on the pixel grid employed by the surrogate model with the functions being constant on each pixel; as a result the relative $L_2$-norm is equivalent to the relative $l_2$-norm.

The dataset we use consists of $\approx 5\,000$ geometries with randomly generated obstacles with $3$, $4$, $5$, $6$, or $12$ edges, with $\approx 1\,000$ geometries for each number of edges.
All computations were performed on NVIDIA V100-GPUs with CUDA 10.1 using Python 3.6 and tensorflow-gpu 2.7~\cite{tensorflow2015-whitepaper}. 

\subsection{Data-based approach} \label{subsec:supervised_approach}

\begin{table}[t]
	\centering
	\begin{tabular}{rr|rrrrr}
		training & \multirow{2}{*}{error} & \multirow{2}{*}{$\frac{\Vert u_{N\!N} - u\Vert_2}{\Vert u\Vert_2}$}   & \multirow{2}{*}{$\frac{\Vert p_{N\!N} - p\Vert_2}{\Vert p\Vert_2}$}  & divergence & momentum & \# \\ 
		data & & &  & residual & residual & epochs \\
		\hline 
		\hline
		\multirow{2}{*}{$10\%$} & training & $2.07\%$   & $10.98\%$ & $1.1\cdot 10^{-1}$ & $1.4\cdot 10^{0}$ & \multirow{2}{*}{$500$} \\
		& validation & $4.48\%$ & $15.20\%$ & $1.6\cdot 10^{-1}$ & $1.7\cdot 10^{0}$ &  \\ \hline
		\multirow{2}{*}{$25\%$} & training & $1.93\%$ & $8.45\%$ & $9.1\cdot 10^{-2}$ & $1.2\cdot 10^{0}$ & \multirow{2}{*}{$500$} \\
		& validation & $3.49\%$ & $10.70\%$ & $1.2\cdot 10^{-1}$ & $1.4\cdot 10^{0}$ &  \\ \hline
		\multirow{2}{*}{$50\%$} & training & $1.48\%$ & $8.75\%$ & $9.0\cdot 10^{-2}$ & $1.1\cdot 10^{0}$ & \multirow{2}{*}{$500$} \\
		& validation & $2.70\%$ & $10.09\%$ &  $1.1 \cdot 10^{-1}$ & $1.2\cdot 10^{0}$ & \\ \hline
		\multirow{2}{*}{$75\%$} & training & $1.43\%$ & $7.30\%$ & $1.0\cdot 10^{-1}$ & $1.5\cdot 10^{0}$ & \multirow{2}{*}{$500$} \\   
		& validation &  $2.52\%$ & $8.67\%$ & $1.2\cdot 10^{-1}$ & $1.5\cdot 10^{0}$ & \\
	\end{tabular}
	\caption{
		\label{tab:supervised_multiple} Performance of the data-based approach  on multiple geometries from the channel dataset compared to OpenFOAM simulations on \textit{locally refined} meshes. The divergence and momentum residuals are averaged over all configurations and pixels.
		}
\end{table}

\begin{figure}[t]
	\centering
	\begin{subfigure}[c]{\textwidth}
		\includegraphics[width=\textwidth]{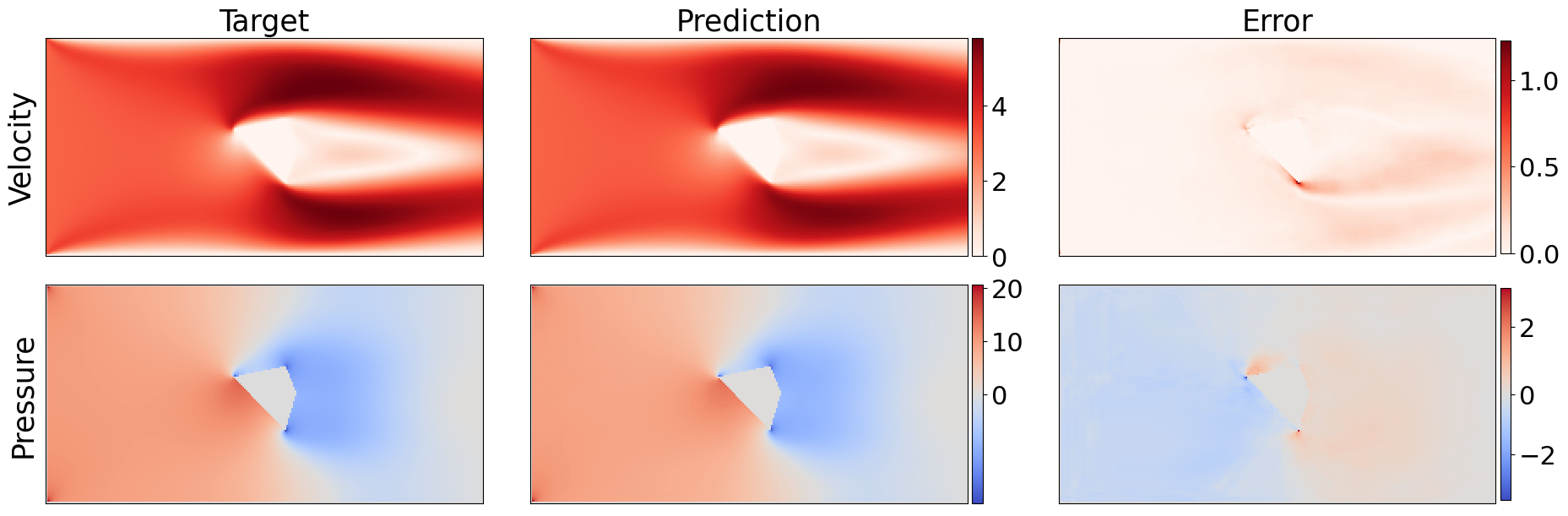} 
		\subcaption{\label{fig:supervised_tp75_geometry_1711} A good prediction. The relative $L_2$-error in $u$ is $2.5\%$ and $6.4\%$ in $p$.
		}
	\end{subfigure}
	\begin{subfigure}[c]{\textwidth}
		\includegraphics[width=\textwidth]{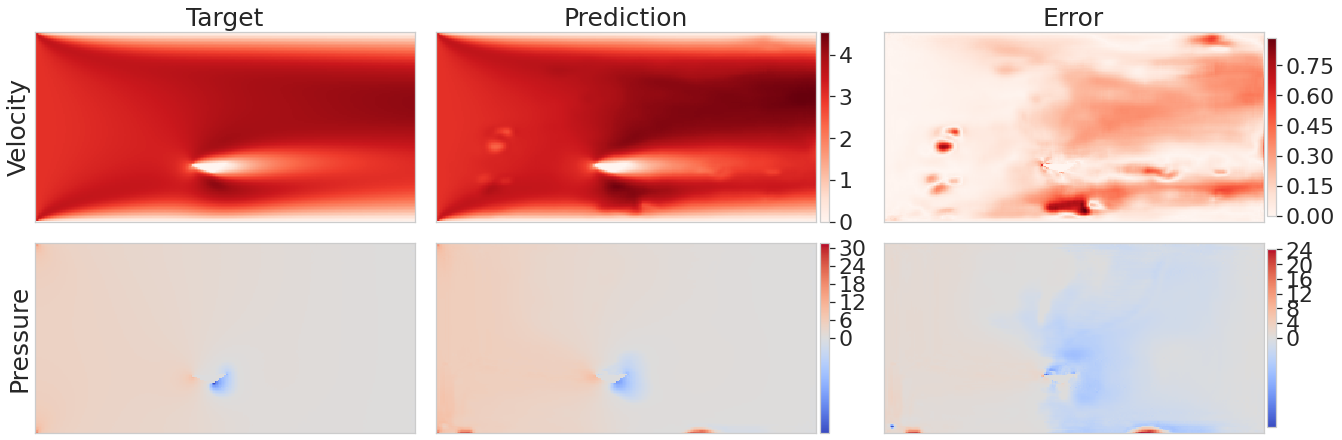} 
		\subcaption{\label{fig:supervised_tp75_geometry_913} A bad prediction. The relative $L_2$-error in $u$ is $5.7\%$ and $28.6\%$ in $p$.
		}
	\end{subfigure}
	\caption{\label{fig:supervised_tp75_predictions} Prediction of velocity and pressure for the data-based approach (Prediction) compared to the OpenFOAM simulation on the \textit{locally refined} mesh (Target). This model was trained on $3750$ geometries. Both geometries are validation geometries.}
\end{figure} 
First, we analyze the fully data-driven approach as discussed in~\cref{sec:cnn}; cf.~\cite{EHK:2021:FPN}. In addition to the velocity field, which was the the focus of~\cite{EHK:2021:FPN}, we also learn the pressure field here; we also did some changes in the network architecture to improve upon the results in~\cite{EHK:2021:FPN}. Later, in~\cref{subsec:unsupervised_approach,subsec:mixed_approach}, we will use the results for the fully data-driven model as a baseline for comparison.

In order to investigate the performance of the data-driven approach, we trained several models with increasing percentages of training data on the channel dataset ($5\,000$ configurations); the remaining data is used for validation, respectively. The training and validation performance, in terms of the relative $L_2$-errors and averaged residual norms, for this approach are summarized in~\cref{tab:supervised_multiple}. We observe that the data-driven model is able to learn the velocity and pressure fields very well, where the errors on the velocity are generally lower. Moreover, the training accuracy is always a bit lower compared to the validation accuracy, which indicates some overfitting. We can observe that using a larger share of training data improves the performance of the model and slightly reduces the overfitting.

For the best model, with $75\,\%$ training data, we also present plots of the velocity and pressure fields  in~\Cref{fig:supervised_tp75_geometry_1711,fig:supervised_tp75_geometry_913}, comparing the reference and prediction data. \Cref{fig:supervised_tp75_geometry_1711} shows a typical example a quantitatively and qualitatively good prediction, whereas the prediction in~\Cref{fig:supervised_tp75_geometry_913} exhibits some clearly unphysical artifacts in the pressure and velocity fields despite a feasible average error. We conjecture that this is due to the pure data loss, which does not include any physical knowledge; as we will discuss in~\cref{sec:results:comparison:data-based}, the physics-aware loss improves this model behavior.

\subsection{Physics-aware approach} \label{subsec:unsupervised_approach}

\begin{table}[t]
	\centering
	\begin{tabular}{rr|rrrrr}
		training & \multirow{2}{*}{error} & \multirow{2}{*}{$\frac{\Vert u_{N\!N} - u\Vert_2}{\Vert u\Vert_2}$}   & \multirow{2}{*}{$\frac{\Vert p_{N\!N} - p\Vert_2}{\Vert p\Vert_2}$}  & divergence & momentum & \# \\ 
		data & & &  & residual & residual & epochs \\ 
		\hline 
		\hline
		\multirow{2}{*}{$10\%$} & training & $4.34\%$ & $9.75\%$ & $2.8 \cdot 10^{-02}$ & $7.4 \cdot 10^{-02}$ & \multirow{2}{*}{$2\,500$} \\
		& validation & $5.70\%$ & $12.81\%$ & $5.7 \cdot 10^{-02}$ & $2.0 \cdot 10^{-01}$ & \\ \hline
		\multirow{2}{*}{$25\%$} & training &$4.17\%$ & $9.61\%$ & $2.5 \cdot 10^{-02}$ & $6.1 \cdot 10^{-02}$ & \multirow{2}{*}{$2\,500$} \\
		& validation & $4.82\%$ & $10.73\%$ & $4.4 \cdot 10^{-02}$ & $1.3 \cdot 10^{-01}$ & \\ \hline 
		\multirow{2}{*}{$50\%$} & training & $4.16\%$ & $9.47\%$ & $2.4 \cdot 10^{-02}$ & $5.7 \cdot 10^{-02}$&  \multirow{2}{*}{$2\,500$} \\
		& validation & $4.37\%$ & $9.68\%$ & $3.7 \cdot 10^{-02}$ & $1.0 \cdot 10^{-01}$ & \\ \hline
		\multirow{2}{*}{$75\%$} & training & $3.82\%$ & $8.71\%$ & $1.8 \cdot 10^{-02}$ & $4.0 \cdot 10^{-02}$ & \multirow{2}{*}{$2\,500$} \\
		& validation & $3.91\%$ & $8.65\%$ & $2.8 \cdot 10^{-02}$ & $8.0 \cdot 10^{-02}$ & \\
	\end{tabular}
	\caption{\label{tab:unsupervised_multiple}Performance of the physics-aware approach on multiple geometries from the channel dataset compared to OpenFOAM simulations on \textit{locally refined} meshes.
		The divergence and momentum residuals are averaged over all configurations and pixels.
	}
\end{table}

\begin{figure}[t]
	\centering
	\begin{subfigure}[c]{\textwidth}
		\includegraphics[width=\textwidth]{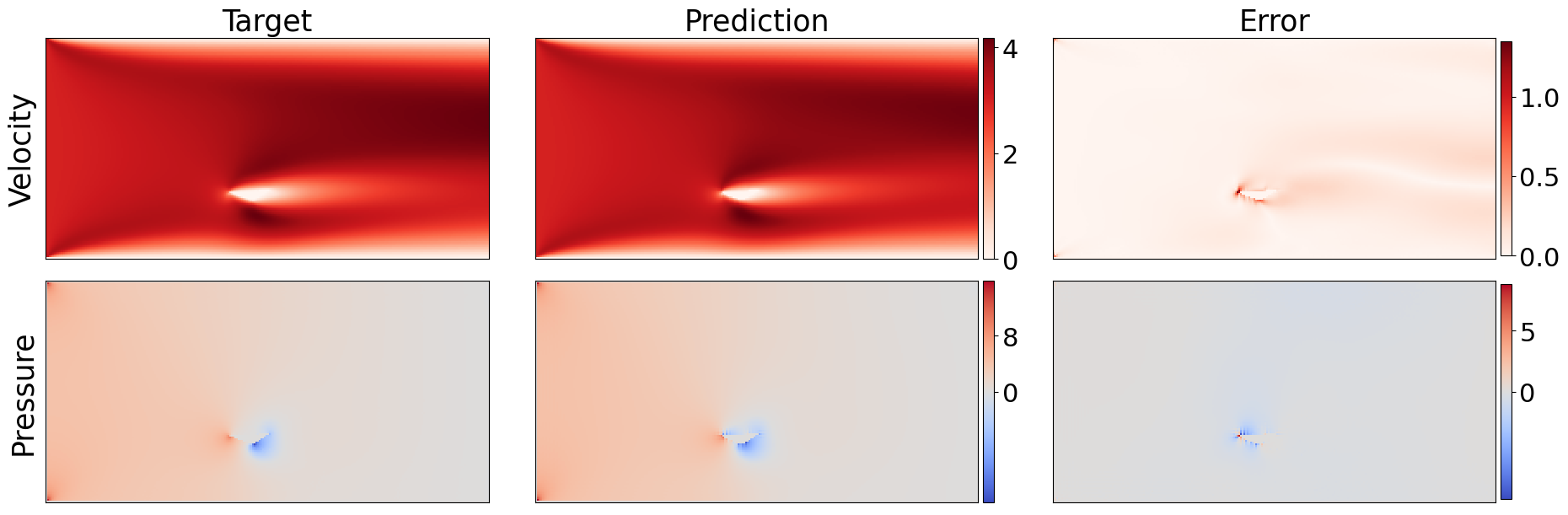} 
		\subcaption{\label{fig:unsupervised_tp75_geometry_913}
			A good prediction. The relative $L_2$-error in $u$ is $2.3\%$ and $6.0\%$ in $p$. 
		} 
	\end{subfigure}
	\vspace{-0.3cm}
	\begin{subfigure}[c]{\textwidth}
		\includegraphics[width=\textwidth]{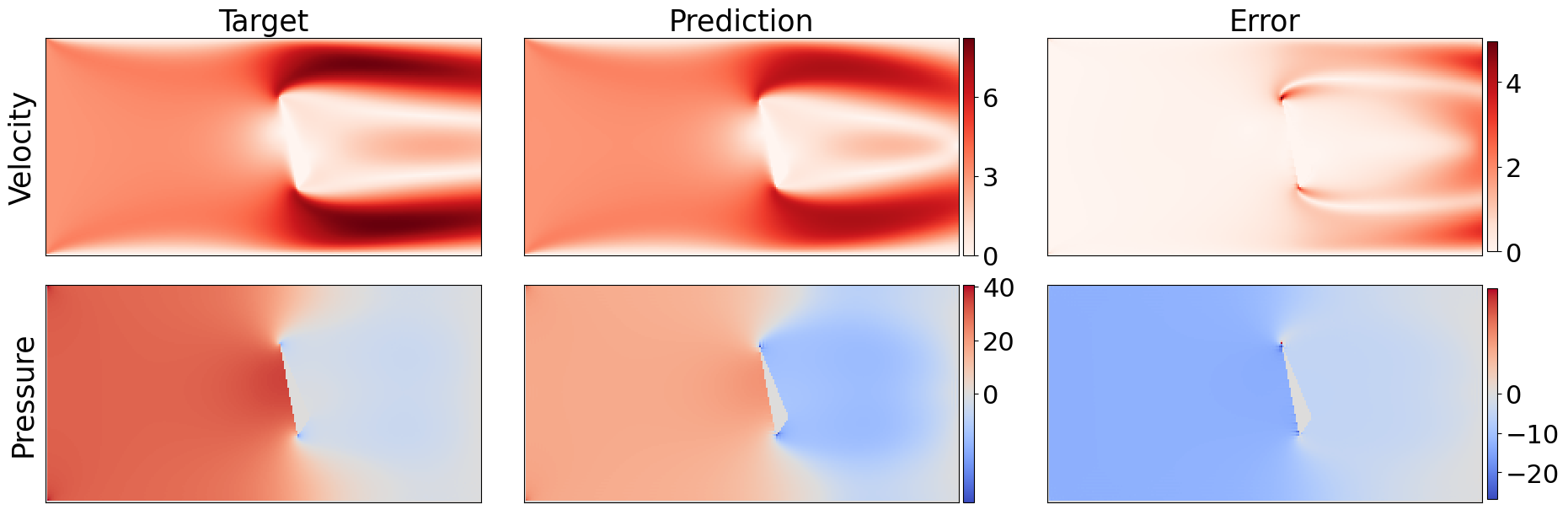} 
		\subcaption{\label{fig:unsupervised_tp75_geometry_1130} The worst prediction. The relative $L_2$-error in $u$ is $21.0\%$ and $45.8\%$ in $p$. 
		}
	\end{subfigure}
	\caption{\label{fig:unsupervised_tp75_predictions}Velocity and pressure for the physics-aware approach (Prediction) compared to the OpenFOAM simulation on \textit{locally refined} meshes (Target). The model was trained on $3\,750$ geometries. All shown geometries are validation geometries.}
\end{figure} 

In this subsection, we will analyze the proposed physics-aware approach in detail. The results on the whole channel dataset ($5\,000$ configurations) for  varying percentages of training data are summarized in~\cref{tab:unsupervised_multiple}. We observe that the physics-aware surrogate model can be extended from a single geometry (\cref{subsec:single-geometry}) to multiple geometries, as discussed theoretically in~\cref{sec:method}. In particular, we obtain predictions with low errors on the velocity and pressure fields, and the performance improves slightly when increasing the share of training data. Interestingly, the overfitting effect is rather small, even when using only $10\,\%$ of the data for training the model.

In \cref{fig:unsupervised_tp75_predictions}, we showcase different predictions from this model. We observe smooth solutions in all cases, without any unphysical artifacts visible. However, when inspecting the error plots, we observe that the error is highest in the vicinity of the obstacle, indicating that the uniform pixel grid cannot fully resolve the boundary layers of the flow. Hence, we observe some error compared with the reference data, which has been computed on a locally refined mesh; cf.~\cref{fig:showcase_mesh}.

In the following, we will discuss the results in more detail: in~\cref{sec:results:comparison:data-based}, we compare the results with the results for the data-based approach in~\cref{subsec:supervised_approach};  in~\cref{sec:high_velocities}, we discuss the correlation of high maximum velocities in the flow field and high prediction errors; and in~\cref{sec:rasterized}, we discuss the influence of the pixel grid on the prediction performance.

\subsubsection{Comparison to the data-based approach} \label{sec:results:comparison:data-based}
\begin{figure}[t]
	\centering
	\begin{subfigure}[c]{0.49\textwidth}
		\includegraphics[width=\textwidth]{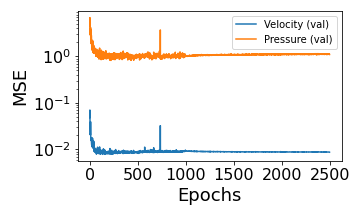} 
		\subcaption{\label{fig:loss_curves_data}Data-based approach.
		} 
	\end{subfigure}
	\begin{subfigure}[c]{0.49\textwidth}
		\includegraphics[width=\textwidth]{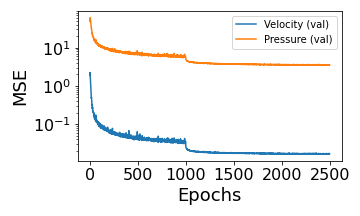} 
		\subcaption{\label{fig:loss_curves_physics}Physics-aware approach
		}
	\end{subfigure}
	\caption{\label{fig:comparison_loss_curves}Validation loss curves for the velocity and pressure over the trained epochs.}
\end{figure} 

Comparing the results in~\cref{tab:supervised_multiple,tab:unsupervised_multiple} for the data-based and the physics-aware approach, respectively, we observe that the data-based model generally has a better performance based on the relative errors compared with the reference data. On the other hand, the residuals of the divergence and momentum equations on the pixel grid are lower for the physics-aware model. This can be easily explained by the fact that the data-based model is trained against the reference data, whereas the physics-aware model is trained to minimize the residuals. 

At first sight, this may seem contradictory since we would expect that, for the same boundary value problem, a lower residual might also results in a lower error. However, as discussed in~\cref{sec:reference_ouput_data}, the reference data is generated based on solving the Navier--Stokes equations with FVM on a locally refined mesh, whereas we evaluate the residuals for the physics-aware model on a uniform pixel grid. This means that the physics-aware model can never reach relative velocity and pressure errors of zero, with our current setting. Likewise, the data-based model will not minimize the residuals on the uniform pixel grid.

Interestingly, the physics-aware model is less prone to overfitting than the data-based model. In particular, despite slightly worse prediction errors, the gap between the training and validation errors is clearly lower for the physics-aware model, indicating better generalization capabilities.

\Cref{tab:supervised_multiple,tab:unsupervised_multiple} also indicate that we performed a significantly larger number of epochs to train the physics-aware than the data-based model on the same data set; in particular, we ran the training for $2\,500$ instead of $500$ epochs. In order to illustrate this, we present plots of the evolution of the mean squared errors for the velocity and pressure over the training process for the physics-aware model with $75\,\%$ training data and a new data-based model that we also trained for $2\,500$ epochs on $75\,\%$ training data in~\cref{fig:comparison_loss_curves}.
Note that this data-based model is not the same model for which we have presented results in this section so far.
The validation errors of the data-based model reach their minimum very quickly, see \cref{fig:loss_curves_data}.
It can be clearly seen that with a training of more than $500$ epochs, the validation errors do not decrease further. On the contrary, they even increase slightly.
Thus, longer training of the data-based model would only lead to stronger overfitting.
In contrast, the validation errors for the physics-aware model decrease more slowly and reach their lowest value only in the further course of the training, see \cref{fig:loss_curves_physics}.

Finally, we briefly discuss those cases where the data-based or the physics-aware model performs badly. As mentioned in~\cref{subsec:supervised_approach}, the data-based approach occasionally makes predictions with unphysical artifacts in the flow and pressure fields. In particular, we show present in~ \cref{fig:supervised_tp75_geometry_913} one example from the validation data set where this is apparent. The corresponding prediction of the physics-aware model for the same sample is shown in~\cref{fig:unsupervised_tp75_geometry_913}. We do not observe the same artifacts. In alignment with the lower overfitting of the physics-aware model, we conclude that the physics-aware model indeed learns better the actual flow behavior based on the residuals of the Navier--Stokes equations. On the other hand, we often see larger errors in the vicinity of the obstacle for the physics-aware approach. This might be attributed due to the uniform pixel grid, which is not specifically refined for resolving the boundary layers in the physics-aware approach. In particular, it seems that the error originates at the obstacle and propagates downstream.

\begin{figure}[t]
	\centering
	\begin{subfigure}[t]{0.48\textwidth}
		\includegraphics[width=\textwidth]{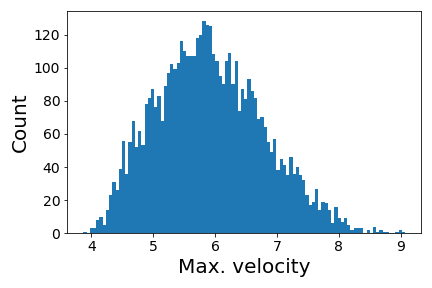}
		\caption{\label{fig:histogramm_max_velo} Histogram of the maximum velocity for the channel dataset}
	\end{subfigure}
	\hfill
	\begin{subfigure}[t]{0.48\textwidth}
		\includegraphics[width=\textwidth]{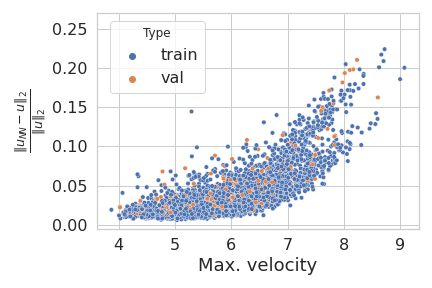} 
		\subcaption{\label{fig:unsupervised_tp75_rel_error} Relative $L_2$-error for $u$  for the physics-aware approach compared to OpenFOAM simulations on \textit{locally refined} meshes. The color indicates training and validation data}
	\end{subfigure}
	
	\vspace{-0.4cm}	
	
	\begin{subfigure}[t]{0.48\textwidth}
		\includegraphics[width=\textwidth]{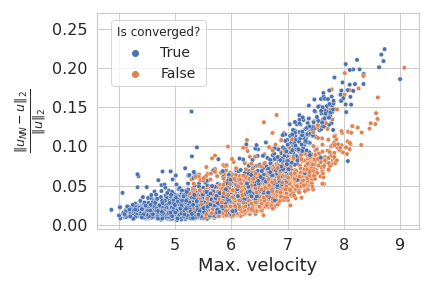} 
		\subcaption{\label{fig:unsupervised_tp75_converged} Relative $L_2$-error for $u$  for the physics-aware approach compared to OpenFOAM simulations on \textit{locally refined} meshes. The color indicates convergence of OpenFOAM simulations on \textit{rasterized} meshes} 
	\end{subfigure}
	\hfill
	\begin{subfigure}[t]{0.48\textwidth}
		\includegraphics[width=\textwidth]{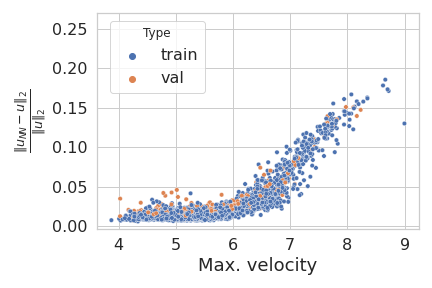} 
		\subcaption{\label{fig:unsupervised_tp75_rasterized_rel_error} Relative $L_2$-error for $u$ for the physics-aware approach compared to OpenFOAM simulations on \textit{rasterized} meshes. The color indicates training and validation data}
	\end{subfigure}
	\caption{Results investigating the correlation of the error with the maximum velocity appearing in the flow field.
		\label{fig:max_vel}}
\end{figure} 

\subsubsection{Correlation of errors and velocities} \label{sec:high_velocities}
In further analyzing the prediction errors, we observed a systematic correlation between the relative error in the velocity and the maximum velocity appearing in the flow field. In particular, depending on the size and position of the obstacle, the maximum velocity can vary significantly; see, e.g., the examples in~\cref{fig:supervised_tp75_predictions,fig:unsupervised_tp75_predictions}. We observe that geometries with a maximum velocity above $6$ exhibit higher average errors compared to those below $6$: $5.8\%$ for $u$ and $12.7\%$ for $p$ versus $2.5\%$ for $u$ and $5.7\%$ for $p$. 

In~\cref{fig:histogramm_max_velo}, we observe that the maximum velocity ranges roughly from $4$ to $9$, following almost a normal distribution; hence, the lower and higher maximum velocities do not appear as often as maximum velocities of $6$. Despite fewer cases with lower maximum velocities, we observe that the prediction error generally increases with an increasing maximum velocity; cf.~\cref{fig:unsupervised_tp75_rel_error}. Moreover, there seems to be no relation to whether the configuration is in the training or validation set. 

Besides arguing based on the distribution of maximum velocities in the data set, it is not surprising that higher maximum velocities lead to higher errors since this might correspond to higher Reynolds numbers and more complex flow patterns. Moreover, our physics-aware loss, as defined in~\cref{eq:ns_pde_loss}, incorporates second-order central stencils for all terms, including the convective terms. 
However, central stencil approximations for convective terms can be problematic when the cell Reynolds number exceeds 2; see, for instance,~\cite[Sec.~2.3]{MCD:2007:LecCFD}. For our uniform pixel grid, this occurs when $\vert u \vert > 4.25 \frac{m}{s}$ in our case. Therefore, it may be necessary to consider alternative approximations for the convective terms. However, this is beyond the scope of this article.

\subsubsection{Influence of the pixel grid} \label{sec:rasterized}
As mentioned before, there seems to be an effect from an insufficient resolution of the boundary layers around the obstacle. In order to investigate potential effects of the resolution of the pixel grid, we rerun all configurations in our data set on the pixel grid; due to their structure, we also denoted these as \textit{rasterized meshes} in~\cref{subsec:single-geometry}.

First of all, we observe that a significant number of OpenFOAM simulations on the rasterized meshes did not converge; cf.~\cref{fig:unsupervised_tp75_converged}. Furthermore, we also see a correlation maximum velocity and convergence in this case. For geometries with obstacles narrower than $1\,\text{m}$ and flow fields with maximum velocities below $6$ $\frac{m}{s}$ almost all simulations converged.
Conversely, for larger obstacles and faster flow fields only a part converged. This is in alignment with our observation on higher errors for higher maximum velocity cases.

Finally, we evaluate the physics-aware model only on those cases where the simulations on the rasterized meshes successfully converged. \Cref{fig:unsupervised_tp75_rasterized_rel_error} displays the errors of the predictions of the physics-aware model against the rasterized simulations. Comparing~\cref{fig:unsupervised_tp75_rel_error,fig:unsupervised_tp75_rasterized_rel_error}, we can observe a much better match when using the rasterized simulations as the reference. For geometries with maximum velocities below $6\,\frac{\text{m}}{\text{s}}$, the average $L_2$-error in $u$ decreases from $2.2\%$ to $1.5\%$. Similarly, for geometries with maximum velocities above $6\,\frac{\text{m}}{\text{s}}$, the average $L_2$ error in $u$ decreases from $6.7\%$ to $5.4\%$. This shows that the pixel grid has an influence on the prediction performance. Further investigations of this aspect are out of the scope of this paper but will be subject of future work.

\subsection{Hybrid approach} \label{subsec:mixed_approach}

\begin{figure}[t]
	\centering
	\includegraphics[width=\textwidth]{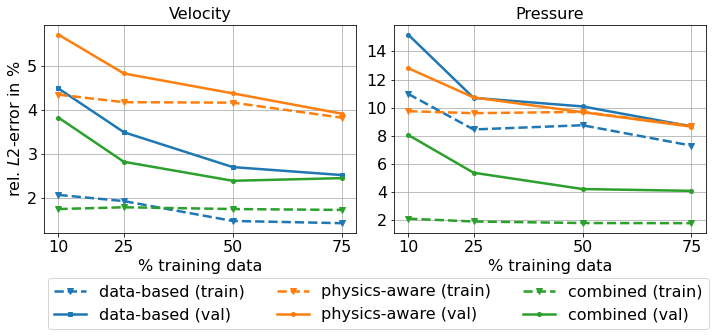}
	\caption{\label{fig:mixed_results_1} Performance of the hybrid approach for abundant simulation results.}
\end{figure}

\begin{figure}[t]
	\centering
	\begin{subfigure}[c]{0.32\textwidth}
		\includegraphics[width=0.9\textwidth]{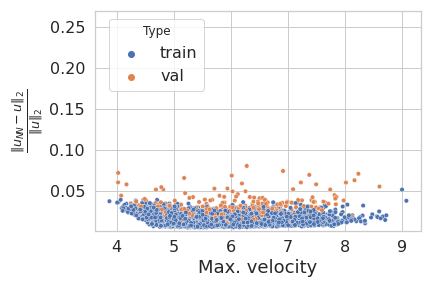} 
		\subcaption{\label{fig:2d_channel_combined_multiple_max_velo_vs_error_summary_data_velo}Data-based velocity.} 
	\end{subfigure}
	\hfill
	\begin{subfigure}[c]{0.32\textwidth}
		\includegraphics[width=\textwidth]{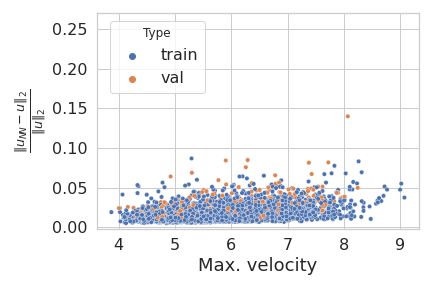} 
		\subcaption{\label{fig:2d_channel_combined_multiple_max_velo_vs_error_summary_physics_velo}Combined velocity.} 
	\end{subfigure}
	\hfill
	\begin{subfigure}[c]{0.32\textwidth}
		\includegraphics[width=\textwidth]{Pictures/Results/unsupervised/multiple/tp75/max_velo_vs_rel_error_u} 
		\subcaption{\label{fig:2d_channel_combined_multiple_max_velo_vs_error_summary_combined_velo}Physics-aware velocity.} 
	\end{subfigure}
	
	\begin{subfigure}[c]{0.32\textwidth}
		\includegraphics[width=\textwidth]{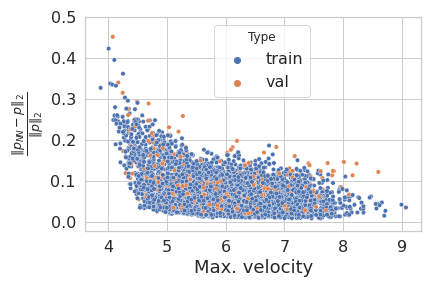} 
		\subcaption{\label{fig:2d_channel_combined_multiple_max_velo_vs_error_summary_data_pressure}Data-based pressure.} 
	\end{subfigure}
	\hfill
	\begin{subfigure}[c]{0.32\textwidth}
		\includegraphics[width=\textwidth]{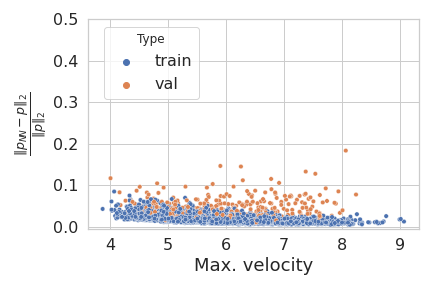} 
		\subcaption{\label{fig:2d_channel_combined_multiple_max_velo_vs_error_summary_physics_pressure}Combined pressure.} 
	\end{subfigure}
	\hfill
	\begin{subfigure}[c]{0.32\textwidth}
		\includegraphics[width=\textwidth]{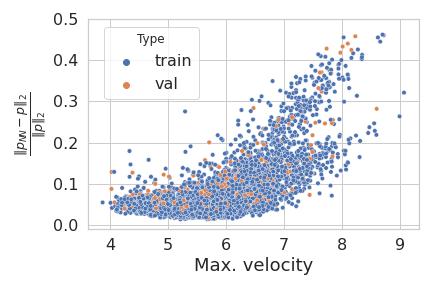} 
		\subcaption{\label{fig:2d_channel_combined_multiple_max_velo_vs_error_summary_combined_pressure}Physics-aware pressure.} 
	\end{subfigure}
	\caption{\label{fig:2d_channel_combined_multiple_max_velo_vs_error_summary}Comparison of the relative $L_2$-error distribution for $u$ and $p$ with regards to the maximum occurring velocity for the data-based ((a) and (d)), combined ((a) and (e)) and physics-aware ((c) and (f)) approaches compared to OpenFOAM simulations on \textit{locally refined} meshes. All models were trained on $3\,750$ geometries.}
\end{figure} 

As discussed in~\cref{sec:results:comparison:data-based}, both approaches have their advantages due to the different loss functions considered. The main advantages of the physics-aware approach are its generalization properties as well as the fact that no reference data is required. The main disadvantage in the CNN approach is that a uniform grid is used, which, in our setting, is not fine enough to fully resolve boundary effects and high velocities. The data-based approach, on the other hand, is able to better capture these. This is presumably because the reference data in the data loss encodes effects which cannot be fully resolved by the pixel grid. However, the data-based model is more prone to overfitting and unphysical flow artifacts.

In order to combine some of the strength of both approaches, we propose a hybrid approach, which employs a weighted sum of the data-based loss function~\cref{eq:mse} and the physics-aware loss function~\cref{eq:ns_pde_loss}. 
This could, for example, be relevant:
\begin{itemize}
	\item if a sufficient number of data samples is available and the generalization properties or physical consistency of the model should be enhanced or
	\item if only an insufficient number of data samples to cover the range of geometries is available; this could specifically be the case if measurement data is used or the simulations are prohibitively expensive. In this case, the missing data can be replaced by using the physics-aware loss.
\end{itemize}

\Cref{fig:mixed_results_1} compares the overall performance of the data-based, physics-aware, and hybrid approaches. Here, the hybrid approach uses both the data-based loss and the physics-aware loss, both with an equal weight of $1$. It can be observed that, for all ratios of training and validation data, the hybrid model outperforms the data-based and physics-aware models in terms of the relative errors in the velocity and pressure; in particular, the prediction performance on the pressure improves significantly, by roughly $50\,\%$, compared with the other approaches. However, unfortunately, the gap between training and validation performance overfitting is on a similar level as for the data-based model.

The performance of the data-based, the physics-aware, and the hybrid approaches with respect to the maximum velocity is shown in~\cref{fig:2d_channel_combined_multiple_max_velo_vs_error_summary}. Whereas the data-based and the physics-aware models show a correlation between the prediction error and the maximum velocity, the hybrid model seems to be rather robust; interestingly, for the data-based approach, we observe a slight deterioration of the performance in the pressure prediction for lower maximum velocities.

The results indicate that, if high-fidelity reference data is available, a combination of the data-based and physics-aware loss functions yields the best results.

\subsection{Weighting of Loss Terms}
There are some elements whose modification may improve the prediction quality of our approach.
This includes, for example, varying the weight of the loss terms, see \cref{eq:ns_pde_loss}.

The initial prediction of the convolutional neural network (CNN) does not fulfill the divergence-free equation due to the random initialization of its weights. During the training process, the minimization of the sum of squared residuals \cref{eq:ns_pde_loss} is pursued, where equal weights ($\omega_{\text{M}} = \omega_{\text{D}} = 1$) may cause the learned prediction to satisfy the momentum equation more than the mass equation. While a valid solution should satisfy both the mass equation and the momentum equation, the mass equation can be seen as primarily serving as a constraint, limiting the space of valid solutions.
Furthermore, in our approach, we employ a variant of the Navier--Stokes equations, specifically the momentum equation, where the assumption $\nabla \cdot \vec{u} = 0$ is explicitly employed to simplify the derivation, as discussed in \cite{gresho1991incompressible}. Therefore, although our primary interest lies in solving the momentum equations, it may be beneficial to confine the search space to velocity fields that comply with the divergence-free condition.
Due to architectural constraints preventing easy modification of our CNNs to guarantee divergence-free predictions, we endeavor to achieve a similar outcome by augmenting the weight $\omega_{\text{D}}$ of the mass residual loss term in the loss function.

\begin{figure}[t]
	\centering
	\includegraphics[width=0.9\textwidth]{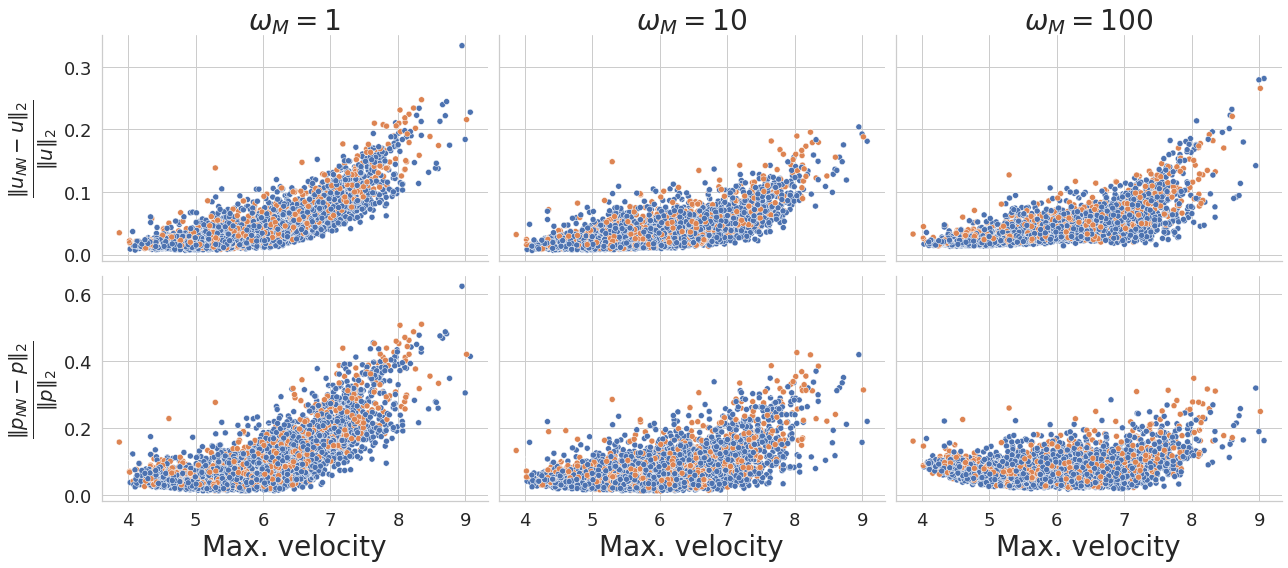} 
	\caption{\label{fig:2d_channel_weight_div}Comparison of the relative $L_2$-error distribution for $u$ and $p$ with regards to the maximum occurring velocity for different approaches compared to OpenFOAM simulations on \textit{locally refined} meshes. All models were trained on $3\,500$ geometries.}
\end{figure} 
Shown in \cref{fig:2d_channel_weight_div} is the distribution of relative errors $L_2$ in velocity and pressure over the maximum occurring velocity for three models for whose training we varied the weight $\omega_{\text{M}}$ of the mass residual in the physics-aware loss from $1$ over $10$ to $100$.
The averaged relative $L_2$ errors are $4.8\%$, $3.7\%$, and $4.5\%$ in the velocity and $10.4\%$, $7.9\%$, and $7.9\%$ in the pressure, for values of $1$, $10$, and $100$ of $\omega_{\text{M}}$, respectively.
Note that we used sixth-order finite difference stencils for all models here, as opposed to second-order stencils in \cref{subsec:unsupervised_approach}.
In addition, the training and validation data sets are not identical to those used in the previous sections. 
Therefore, the error values reported in this section are not necessarily directly comparable to the previous ones.

With an increase in the weight of the mass residual, we see a reduction in the error in the velocity as well as in the pressure, especially at higher velocities.
This effect is very clear for $\omega_{\text{M}} = 10$.
In the averaged errors, this model improves by about $1\%$ in velocity and about $2.5\%$ in pressure compared to the model trained with equal weights, i.e. $\omega_{\text{M}} = 1$.
However, for the model trained with $\omega_{\text{M}} = 100$ we see larger errors in the velocity overall. 
Here, even for low velocities, the errors in the pressure become larger.
An even further increase of the weight $\omega_{\text{M}}$ led to a deterioration of the predictive capabilities, because while the predictions of the model increasingly satisfied the divergence-free constraint, the momentum residual grew.

We have thus demonstrated that increasing the weight of the mass residual in the physics-aware loss can significantly improve the predictive capabilities of the model.

\subsection{Test Data}
So far, our predictions have been limited to geometries that were present in the training or validation datasets. These datasets exclusively pertain to the model problem, as illustrated in \cref{fig:showcase_domain}. Notably, these geometries encompass obstacles in the form of star-shaped polygons with up to 12 vertices.
In this section, we will showcase predictions obtained using the physics-aware convolutional neural network for geometries that possess alternative types of obstacles.
In doing so, we assess how well the model can generalize to previously unseen geometries, extrapolate beyond the training data, and effectively handle new and unique shapes.

\begin{figure}[t]
	\centering
	\begin{subfigure}[c]{\textwidth}
		\includegraphics[width=\textwidth]{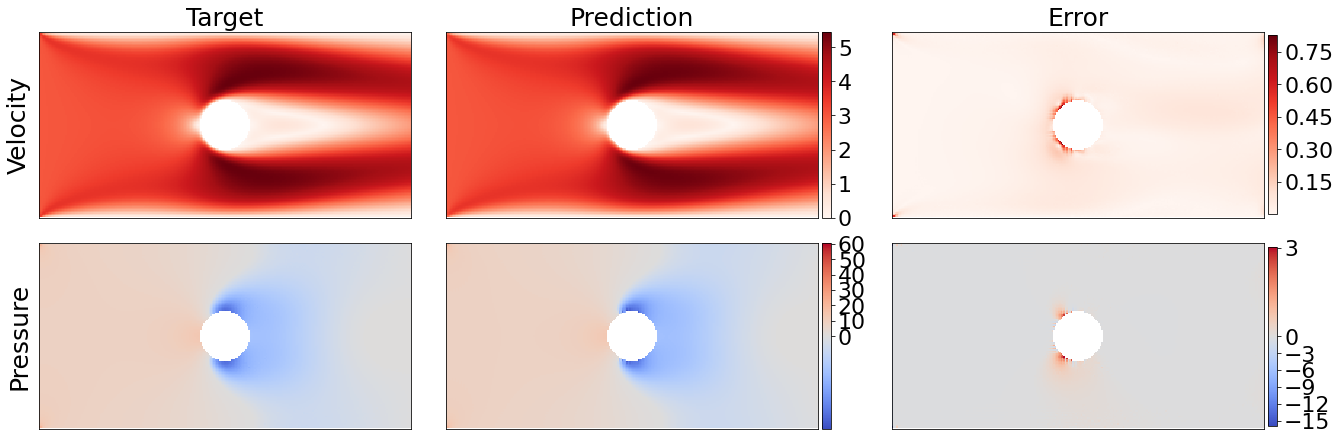}
		\subcaption{\label{fig:unsupervised_tp75_test_prediction_circle_4}
			A circular obstacle. The relative $L_2$-error in $u$ is $1.2\%$ and $4.9\%$ in $p$.
		} 
	\end{subfigure}
	\vspace{-0.3cm}
	\begin{subfigure}[c]{\textwidth}
		\includegraphics[width=\textwidth]{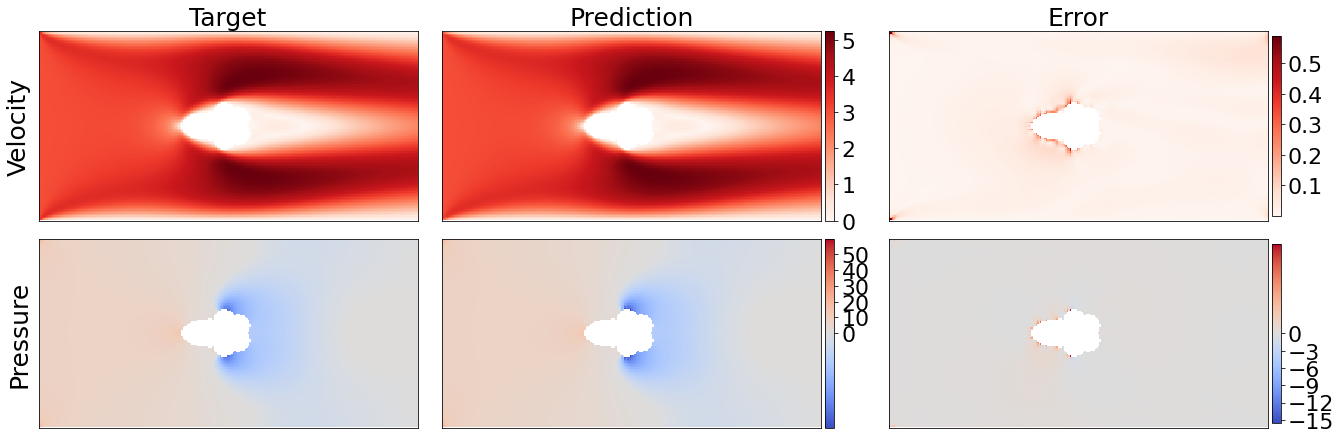}
		\subcaption{\label{fig:unsupervised_tp75_test_prediction_composite_25}
			A composite obstacle. The relative $L_2$-error in $u$ is $0.7\%$ and $4.8\%$ in $p$.
		}
	\end{subfigure}
	\caption{\label{fig:unsupervised_tp75_test_prediction}Velocity and pressure for the physics-aware approach (Prediction) on test geometries compared to the OpenFOAM simulation on \textit{locally refined} meshes (Target). The model was trained on $3\,750$ geometries. All shown geometries are validation geometries.}
\end{figure} 

The first geometry we will test the model on is a circle with radius $0.4$ that we place in the middle of the channel.
This type of obstacle is a highly distinct and different type of obstacle compared to the star-shaped polygons present in the training dataset.
Circles have a continuous curved boundary, which contrasts with the sharp edges of the star-shaped polygons, making them significantly novel geometries for the model.
We show the prediction of our model for this geometry in \cref{fig:unsupervised_tp75_test_prediction_circle_4}. Smooth predictions are observed. High errors occur only near the obstacle. This prediction shows that our model can handle the curvature of a circle very well without having seen a single curved obstacle during training.

The second geometry we test our model on is a composition of an oval and a flower with 5 petals. The oval has a horizontal radius of $0.45$ and a vertical radius of $0.25$. The flower has a maximum radius of $0.4$.
The 5-petaled flower also differs from a circle because the curvature is not uniform throughout, but is interrupted by sharp bends where the petals meet.
We show the prediction of our model for the second geometry in \cref{fig:unsupervised_tp75_test_prediction_composite_25}. Again, we see smooth predictions with high errors occurring only near the obstacle.

These two predictions exemplify that our model is capable of making reasonable and accurate predictions for geometries with significantly different obstacles.

\subsection{Computation Time} \label{sec:computation_time}
	An important aspect of a surrogate model is the speed with which it can be evaluated.
	Therefore, in this section we want to compare the time needed to evaluate the surrogate model with the time needed for a reference CFD simulation.
	A CFD simulation described in \cref{sec:reference_ouput_data} takes between $10$ and $60$ minutes, depending on the geometry, and in individual difficult cases the simulation may take longer than $60$ minutes.
	In comparison, we need only roughly $6$ milliseconds (ms) to evaluate our surrogate models on a geometry. Thus, the evaluation of the surrogate model is between $100\,000$ and $600\,000$ times faster than a CFD simulation.
	This does not include the time required to mesh the geometry, which can be very time consuming depending on the geometry and mesh fineness, and to set up the CFD simulation.
	
	However, the time required for training the surrogate model is very high.
	For example, one training step for one geometry takes about $50$ ms.
	Thus, training the model discussed in \cref{subsec:unsupervised_approach}, which was trained on about $3\,750$ geometries, takes roughly $5$ days.
	The training, however, can be done in an offline phase before the actual deployment of the surrogate model, so the long training time is not as significant, especially when many simulations need to be run quickly.

\section{Conclusion}
\label{sec:conclusion}
We have introduced a novel physics-aware approach to train convolutional neural networks as surrogate models that relies exclusively on the physics modeling fluid behavior in multiple irregular geometries. 
Our approach does not rely on reference data and only requires the geometry image and boundary conditions. 
However, we have demonstrated that incorporating the physics-aware loss in the training process improves upon the data-based approach when reference data is available.
This approach serves as an excellent surrogate model, with the evaluation being in the order of $O\left( 10^5 \right)$ times faster than a conventional CFD simulation.

Our physics-aware approach performs well for low velocity geometries and demonstrates strong generalization capabilities. 
In contrast, the data-based approach struggles to generalize for the same low velocity geometries. 
Despite using a coarser resolution and finite differences, which may not be ideal for Navier--Stokes, our models achieve excellent predictions close to the reference solution for most geometries. 
However, for cases where our physics-aware models did not match the reference solution, even higher resolution finite volume methods failed to obtain a converged solution. 
Consequently, accurate predictions cannot be expected in such cases.

\section*{Acknowledgments}
This work was performed as part of the Helmholtz School for Data Science in Life, Earth and Energy (HDS-LEE) and received funding from the Helmholtz Association of German Research Centers. 
We gratefully acknowledge the use of the computational facilities of the Center for Data and Simulation Science (CDS) at the University of Cologne and of the Department of Mathematics and Computer Science of the Technische Universit\"at Bergakade\-mie Freiberg  operated by the University Computing Center (URZ) and funded under grant application No. $100376434$ to the State Ministry for Higher Education, Research and the Arts (SMWK) of the Federal State of Saxony on Artificial Intelligence and Robotics for GeoEnvironmental Modeling and Monitoring.

\bibliographystyle{siamplain}
\bibliography{GeoPINN}
\end{document}